\documentclass[10pt,final]{siamltex}
\usepackage{amsmath}
\usepackage{bm}
\usepackage{amssymb,version}
\usepackage{cases}
\usepackage{color}
\usepackage{verbatim}

\usepackage{graphicx}
\usepackage{subfig}
\allowdisplaybreaks

\newcommand{\normmm}[1]{{\left\vert\kern-0.25ex\left\vert\kern-0.25ex\left\vert #1
    \right\vert\kern-0.25ex\right\vert\kern-0.25ex\right\vert}}
\newcommand{\ttt}[1]{(-\frac{t^{n+1}}{T})}

\begin{document}
    \title{ A class of new linear, efficient and high-order implicit-explicit methods for the coupled free flow-porous media system based on nonlinear Lions interface condition
\thanks{This work is supported in part by the National Natural Science Foundation of China  grants 12271302, 12131014 and 12271303, and the Major Fundamental Research Project of Shandong Province of China grant ZR2023ZD33, and the Taishan Scholars Program of Shandong Province of China grant tsqn201909044.}}

    \author{ Xinhui Wang\thanks{School of Mathematics, Shandong University, Jinan, Shandong, 250100, P.R. China. Email: xinhuiwang0108@163.com}.
         \and Xu Guo\thanks{Geotechnical and Structural Engineering Research Center, Shandong University, Jinan, Shandong, 250100, P.R. China. Email: guoxu@sdu.edu.cn}.
        \and Xiaoli Li\thanks{Corresponding Author. School of Mathematics, Shandong University, Jinan, Shandong, 250100, P.R. China. Email: xiaolimath@sdu.edu.cn}.
}

\maketitle

\begin{abstract}
In this paper, we construct and analyze new first- and second-order implicit-explicit (IMEX) schemes for the unsteady Navier-Stokes-Darcy model to describe the coupled free flow-porous media system, which is based on the scalar auxiliary variable (SAV) approach in time and finite element method in space. The constructed schemes are linear, only require solving
 a sequence of linear differential equations with constant coefficients at each time step, and can decouple the Navier-Stokes and Darcy systems. The unconditional stability of both the first- and second-order IMEX schemes can be derived for the coupled system equipped with the Lions interface condition, where the key point is that we should construct a new trilinear form to balance the fully explicit discretizations of the nonlinear terms in the complex system. We can also establish rigorous error estimates for the velocity and hydraulic head of  the first-order scheme without any time step restriction. Numerical examples are presented to validate the proposed schemes.
\end{abstract}


\section{Introduction}
The study of fluids flowing coupled porous media systems has been successful in a wide range of industrial and geophysical applications, and has received a lot of attention in the literature \cite{Hill2008Poiseuille,M2019convection,he2015domain}. The Navier-Stokes-Darcy model is of great significance in simulating the coupling of free fluid flow and flow in porous media, such as surface and groundwater flow \cite{dis2004domain,Layton2002Coupling}, protection of karst aquifers \cite{cao2010coupled}, and industrial filtration \cite{Hanspal2006numeri}. The main difficulties that this model has to face are the nonlinear convection term and the coupled interface conditions.

Numerical solutions of free flow-porous media systems play an important role in practical application and an enormous amount of works have been devoted on the design, analysis and implementation of numerical schemes, including steady Navier-Stokes-Darcy model \cite{badea2010numeri,cai2009numeri, Chidyagwai2009on, cao2013decoupling, Girault2009DG, he2015domain,Shiue2018Convergence} and unsteady Navier-Stokes-Darcy model \cite{ce2008analysis,ce2009Primal,ce2013time,Qiu2020domain}. In this paper we focus on the unsteady Navier-Stokes-Darcy model. Rui and Zhang \cite{RUI2017stabilized} constructed a stabilized mixed finite element method for solving the coupled Stokes and Darcy flow equations with a solute transport. They proposed a mixed weak formulation and used the nonconforming piecewise Crouzeix-Raviart finite element, piecewise constant and conforming piecewise linear finite element to approximate velocity, pressure and concentration respectively. Chen et al. \cite{Chen2013Efficient} proposed and studied two second-order in time implicit-explicit methods for the coupled Stokes–Darcy system and established the unconditional and uniform in time stability for both schemes. They also \cite{Chen2020Uniquely} proposed a family of fully decoupled numerical schemes such that the Navier-Stokes equations, the Darcy equations, the heat equation and the Cahn-Hilliard equation are solved independently at each time step for the extended  Cahn-Hilliard-Navier-Stokes-Darcy-Boussinesq system. A domain decomposition method was proposed to solve a time-dependent Navier-Stokes-Darcy model with Beavers-Joseph interface condition and defective boundary condition by Qiu et al. \cite{Qiu2020domain}. The convergences of this domain decomposition method were rigorously analyzed for the time-dependent Navier-Stokes-Darcy model with Beavers-Joseph interface condition.

As far as we know, all above works are based on the implicit or the implicit-explicit discretizations for the nonlinear convection term for the Navier-Stokes equation, thus one should solve a nonlinear system or the linear system with variable coefficients at each time step. From a computational point of view, it would be ideal to be able to treat the nonlinear term explicitly without any stability constraint.  To the best of our knowledge,  the developed IMEX schemes for the Navier-Stokes equation in \cite{lin2019numerical,li2020error}, based on the scalar auxiliary variable (SAV) approach \cite{shen2018scalar,shen2018convergence}, can maintain unconditionally energy stability with explicit treatment of nonlinear term. In our recent work \cite{li2022new,li2022fully}, we construct the new SAV scheme with purely linear and better stability, which can be regarded as an enhanced version of the original SAV method in \cite{lin2019numerical}. Following the SAV ideas in \cite{li2022new}, Jiang and Yang \cite{jiang2021SAV,jiang2023fast} proposed and studied efficient artificial compressibility ensemble schemes for fast computation of Stokes-Darcy flow ensembles. To the best of our knowledge, there is no related work for the Navier-Stokes-Darcy model with explicit discretization for the nonlinear convective and coupled interface terms to still maintain unconditional stability, much less convergence analysis.

So the aim of this work is to extend the approach proposed in \cite{li2022new,jiang2021SAV} to the unsteady Navier-Stokes-Darcy model, which is much more complicated with nonlinear convective and coupled interface terms between the free flow and porous media systems, and establish rigorous error estimates for the constructed schemes without any time step restriction. The main contributions of this paper are shown in:
\begin{itemize}
\item[$\bullet$] We construct first- and second-order implicit-explicit (IMEX) schemes for the unsteady Navier-Stokes-Darcy model, where the constructed schemes are linear, only require solving
 a sequence of linear differential equations with constant coefficients at each time step, and can decouple the Navier-Stokes and Darcy systems.

\item[$\bullet$] Inspired by the new SAV approach in \cite{li2022new}, unconditional stability of both first- and second-order IMEX schemes can be derived, where the key procedure is based on the fact that nonlinear interface terms need to be treated more carefully and we should construct a new trilinear form to balance the fully explicit discretizations of the nonlinear terms
for the highly coupled system equipped with the Lions interface condition.

\item[$\bullet$] By using the built-in stability result, we can establish rigorous error estimates for the velocity and hydraulic head of the first-order scheme without any time step restriction. Note that the error analysis here is much more difficult than that in \cite{li2022new} due to its highly coupled and nonlinear nature.

\end{itemize}

The framework of this paper is as follows. The next section describes the unsteady Navier-Stokes-Darcy model and gives its weak form. In section 3, we construct the first- and second-order schemes and derive the unconditional stability. In Section 4, we establish the rigorous error estimates for the constructed first-order scheme. Subsequently, we support our theoretical analyses with some numerical examples to confirm the efficiency of the algorithms in Section 5. Finally, we summarise the main points of the paper.

\section{The unsteady Navier-Stokes-Darcy model}
In this section, we first introduce the unsteady Navier-Stokes-Darcy model, including the governing equations, interface conditions and weak forms.

\begin{figure}
  \centering
  \includegraphics[width=7.1cm]{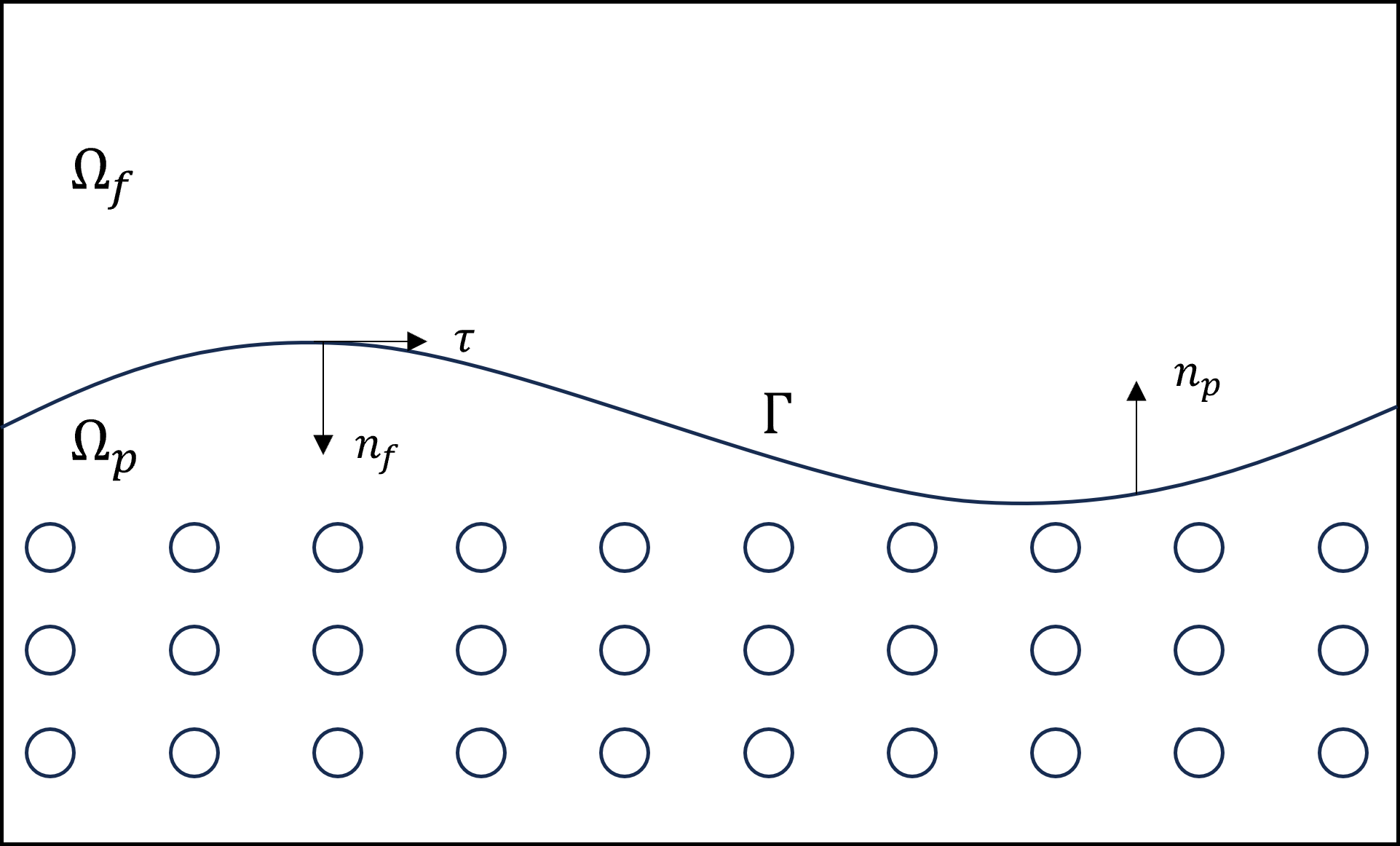}
  \caption{Coupled domains of free flow and porous media flow}\label{NSDdomain}
\end{figure}

We consider the bounded domain $\Omega\subset \mathbb{R}^d~(d=2,3)$ which is divided into two disjoint domains $\Omega_f$ and $\Omega_p$. Let $\Gamma$ denote the interface between two domains $\Gamma=\partial\Omega_f\cap\partial\Omega_p$, see Figure \ref{NSDdomain}.
The fluid flow in domain $\Omega_f$ is governed by the Navier-Stokes equations
\begin{align}
u_t-\nu\Delta u+(u\cdot\nabla) u+\nabla p&=f_1, \ \ \  {\rm in}~\Omega_f\times(0,T],\label{stokes1}\\
\nabla\cdot u&=0, \ \ \ \  {\rm in}~\Omega_f\times(0,T],\label{stokes2}\\
u(x,0)&=u_0, \ \ \ {\rm in}~\Omega_f,\label{stokes t=0}\\
u(x,t)&=0, \  \ \ \  {\rm on}~\partial\Omega_f\backslash\Gamma\times (0,T].\label{stokes gamma}
\end{align}
The coefficient $\nu$ is the kinematic viscosity, $T$ is the final time. $u$ and $p$ represent the velocity and the pressure of free flow in $\Omega_f$, respectively. The function $f_1$ is an external body force.

The flow in porous media domain $\Omega_p$ satisfies the Darcy's law
\begin{align}
u_p&=-\mathcal{K}\nabla \phi, \ \ \  {\rm in}~\Omega_p\times(0,T],\label{darcy1}\\
S_0\phi_t+\nabla\cdot u_p&=f_2, \ \ \ \qquad {\rm in}~\Omega_p\times(0,T],\label{darcy2}\\
\phi(x,0)&=\phi_0, \ \ \ \ \ \ \ \ \ {\rm in}~\Omega_p,\label{darcy t=0}\\
\phi(x,t)&=0, \ \ \ \ \ \ \ \ ~~~{\rm on}~\partial\Omega_p\backslash\Gamma\times (0,T].\label{darcy gamma}
\end{align}

$u_p$ represents the velocity of the porous media flow, $\phi$ is the hydraulic head in $\Omega_p$, $S_0$ is a specific mass storativity coefficient and $f_2$ is the source/sink term. $\mathcal{K}=k\mathbb{I}$ is the hydraulic conductivity tensor, which is a symmetric positive definite matrix. (\ref{darcy1})-(\ref{darcy2}) can be rewritten as follows
\begin{align}
S_0\phi_t-\nabla\cdot (\mathcal{K}\nabla \phi)&=f_2, \qquad {\rm in}~\Omega_p\times(0,T].\label{darcy}
\end{align}
From now on, without causing confusion, we will use $u$ to represent $u_f$. The Navier-Stokes and Darcy systems are coupled by the following interface conditions:
\begin{align}
u_f\cdot{n}_f+u_p\cdot{n}_p&=0, \qquad~~{\rm on}~\Gamma\times (0,T],\label{gamma1}\\
p-\nu{n}_f\frac{\partial u}{\partial {n}_f}+\frac{1}{2}u\cdot u&=g\phi, \qquad{\rm on}~\Gamma\times (0,T],\label{gamma2}\\
-\nu{\tau}_i\frac{\partial u}{\partial {n}_f}&=\alpha\sqrt{\frac{\nu g}{tr(\mathcal{K})}}(u\cdot{\tau}_i),~i=1,\cdots,d-1, \qquad{\rm on}~\Gamma\times (0,T],\label{gamma3}
\end{align}
where ${n}_f$ and ${n}_p$ denote the unit outward vectors on $\partial \Omega_f$ and $\partial \Omega_p$ respectively. It is easy to see that ${n}_f=-{n}_p$. $g$ is the gravitational acceleration, ${\tau}_i(i=1,\cdots,d-1)$ represent the unit tangential vectors on $\Gamma$, $\alpha$  is an experimentally derived data related to the properties of the porous media domain.

The first equation is the mass conservation, as well as the velocity continuity across the interface $\Gamma$. The second equation, also called the Lions interface condition, is the balance of the normal stress, the velocity of free fluid is equal to the normal component of the normal velocity in the porous media domain and considering the inertial force \cite{m2022predicting,Chidyagwai2009on,Girault2009DG,M2019convection}.
The inclusion of the $\frac{1}{2}u\cdot u$ term in this interface condition is essential in conducting the nonlinear stability analysis.
The last equation is the well-known Beavers-Joseph-Saffman condition, which indicates that the velocity is proportional to the normal shear stress in the tangential direction.

Let us denote the $L^2$ inner product in $\Omega_i$ by $(\cdot,\cdot)_i,i=f,p$,
and define some Sobolev spaces
\begin{align*}
H_f&=\{v\in H^1({\Omega_f})^d:v|_{\partial\Omega_f\backslash\Gamma}=0\},\\
H_p&=\{\psi \in H^1(\Omega_p):\psi|_{\partial\Omega_p\backslash\Gamma}=0\},\\
Q_f&=L^2(\Omega_f),\\
W&=H_f\times H_p.
\end{align*}


The weak formulation of the unsteady Navier-Stokes-Darcy model is described as follows: $\forall~t\in (0,T]$, find $\bm{u}=(u,\phi)\in W,~p\in Q_f$ such that
\begin{equation}\label{nsd problem}
\left\{
\begin{aligned}
&(u_t,v)_f+gS_0(\phi_t,\psi)_p+a_N(u,u,v)+\nu(\nabla u,\nabla v)_f\\
&\qquad+\alpha\sqrt{\frac{\nu g}{tr(\mathcal{K})}}\sum_{i=1}^{d-1}
\int_\Gamma(u\cdot\tau_i)(v\cdot\tau_i)ds
+g(\mathcal{K}\nabla\phi,\nabla\psi)_p
+c_\Gamma({v},\phi)
-c_\Gamma({u},\psi)\\
&\qquad
-(p,\nabla\cdot v)_f=(f_1,v)_f+g(f_2,\psi)_p,\qquad\forall\bm{v}=(v,\psi)\in W,\\
&(\nabla\cdot u,q)_f=0,\qquad\qquad~\forall q\in Q_f,\\
&(u(0),v)_f=(u_0,v)_f,\qquad\forall v\in H_f,\\
&(\phi(0),\psi)_p=(\phi_0,\psi)_p,\qquad\forall \psi\in H_p,
\end{aligned}
\right.
\end{equation}
where
\begin{align*}
&a_N(u,v,w)=((u\cdot\nabla) v,w)_f- \frac{1}{2} \int_\Gamma ( u\cdot v )  ( w\cdot n_f ) ds, \qquad\forall u,v,w\in H_f,\\
&c_\Gamma(v,\phi)=g\int_\Gamma\phi v\cdot n_f ds, \qquad\forall v\in H_f,~\phi\in H_p.
\end{align*}

We now give some notations:
\begin{align*}
& \ \normmm{v}_{m,k,i}=\left(\Delta t\sum_{n=0}^{N}\|v^n\|^m_{H^k(\Omega_i)}\right)^{1/m},  \ \ \normmm{v}_{\infty,k,i}=\max_{0\leq n\leq N}{\|v^n\|_{H^k(\Omega_i)}}, \\
&\ \ \|v\|_{m,k,i}=\|v\|_{L^m(0,T;H^k(\Omega_i))},
\end{align*}
and $\|\cdot\|_{k,i}=\|\cdot\|_{H^k(\Omega_i)}$, where $v^n=v(t^n)$, $i=f,p$.

\medskip
\begin{lemma}\label{lemma 1}
For the trilinear term $a_N(u,v,w)$, the following inequalities hold \cite{shen1992on,Connors2012fluid}:
\begin{align}
&a_N(u,v,w)\leq C\|u\|_{1,f}\|v\|_{1,f}\|w\|_{1,f},\\
&a_N(u,v,w)\leq C\|u\|_{0,f}\|v\|_{2,f}\|w\|_{1,f},\\
&|\int_\Gamma (u\cdot v) (w\cdot n_f) ds|\leq
C\|u\|_{0,f}^{1/2}\|u\|_{1,f}^{1/2}
\|v\|_{0,f}^{1/2}\|v\|_{1,f}^{1/2}
\|w\|_{0,f}^{1/2}\|w\|_{1,f}^{1/2},\label{gamma term}
\end{align}
and when take $d=2$,
\begin{align}
a_N(u,v,w)\leq C\|u\|_{0,f}^{1/2}\|u\|_{1,f}^{1/2}
\|v\|_{0,f}^{1/2}\|v\|_{1,f}^{1/2}\|w\|_{1,f},
\end{align}
where $C$ is a positive constant which is dependent on $\Omega_f$.
\end{lemma}

\medskip
\begin{lemma}\label{lemma2}
Suppose that $\nabla\cdot u=0$, the trilinear term has the following property
\begin{align} \label{an=0}
a_N(u,u,u)=0,\qquad u\in H_f,~\nabla\cdot u=0.
\end{align}

\begin{proof}
 By using Green's formula, we have
   \begin{align*}
  (\nabla\cdot u,v\cdot w)_f=-(u,\nabla(v\cdot w))_f+\int_\Gamma(v\cdot w)(u\cdot n_f)ds,
  \qquad \forall u,v,w\in H_f,
  \end{align*}
  where
  \begin{align*}
  (u,\nabla(v\cdot w))_f
  =((u\cdot\nabla) v,w)_f+((u\cdot\nabla) w,v)_f.
  \end{align*}
We can therefore deduce that
  \begin{align*}
  ((u\cdot\nabla) v,w)_f=-(\nabla\cdot u,v\cdot w)_f-((u\cdot\nabla) w,v)_f+\int_\Gamma(v\cdot w)(u\cdot n_f)ds,\qquad \forall u,v,w\in H_f.
  \end{align*}
  Since $\nabla \cdot u=0$, we have
    \begin{align*}
  ((u\cdot\nabla) v,w)_f=-((u\cdot\nabla) w,v)_f+\int_\Gamma(v\cdot w)(u\cdot n_f)ds,\qquad \forall u,v,w\in H_f,~\nabla\cdot u=0.
  \end{align*}
  Hence
   \begin{align*}
  ((u\cdot\nabla) v,w)_f = &\frac{1}{2}  ((u\cdot\nabla) v,w)_f + \frac{1}{2}  ((u\cdot\nabla) v,w)_f \\
 =& \frac{1}{2}\left(((u\cdot \nabla )v,w)_f-((u\cdot \nabla)w,v)_f+\int_\Gamma(v\cdot w)(u\cdot n_f)ds\right).
  \end{align*}
  Thus we obtain the desired result \eqref{an=0} by setting $v=w=u$.
\end{proof}

\end{lemma}

We then construct regular triangulations $\tau_h^f$ and $\tau_h^p$ of $\Omega_f$ and $\Omega_p$ respectively, which are compatible in the interface $\Gamma$. Define the following finite element spaces
\begin{align*}
H_{fh}&=\{v\in H_f: v|_e\in (P_k)^d,~\forall e\in \tau_h^f\},\\
H_{ph}&=\{\psi\in H_p : \psi|_e\in P_m,~\forall e\in \tau_h^p\},\\
Q_{fh}&=\{q\in Q_f: q|_e\in P_j,~\forall e\in \tau_h^f\},\\
W_h&=H_{fh}\times H_{ph}.
\end{align*}
Also, we define the divergence-free space
\begin{align*}
V_h=\{v_h\in H_{fh}:(q_h,\nabla\cdot v_h)_f=0,\forall q_h\in Q_{fh}\}.
\end{align*}

\section{Numerical schemes}
In this section, we construct first- and second-order linear and efficient schemes  for the unsteady Navier-Stokes-Darcy model  and derive the unconditional stability for the constructed schemes.

At the beginning, inspired by \cite{li2022new,li2022fully}, we introduce an auxiliary variable $r(t^{n+1})=exp(-\frac{t^{n+1}}{T})$, which is only relevant to the time direction. Then we can derive that
\begin{align*}
r_t=-\frac{1}{T}r+\frac
{1}{exp\ttt)}(c_\Gamma(u,\phi)-c_\Gamma(u,\phi)
+a_N(u,u,u)).
\end{align*}
The last three terms of the above equation is equal to 0 by using Lemma \ref{lemma2}. In the discrete format, however, the presence of these terms plays an important role in ensuring the stability of the numerical schemes and in decoupling the Navier-Stokes and Darcy systems.
For the sake of simplicity, we define
\begin{align*}
\eta=\alpha\sqrt{\frac{\nu g}{tr(\mathcal{K})}}.
\end{align*}

Assume that the solution of problem (\ref{nsd problem}) has these regularity results:
\begin{equation}\label{regularity solution}
\begin{cases}
u\in L^\infty(0,T;H^{k+1}(\Omega_f)),\quad \partial_t u\in L^2(0,T;H^{k+1}(\Omega_f)),\quad
\partial _{tt} u\in L^2(0,T;L^2(\Omega_f)),\\
\phi\in L^\infty(0,T;H^{m+1}(\Omega_p)),\quad \partial_t\phi\in L^2(0,T;H^{m+1}(\Omega_p)),
\quad \partial_{tt}\phi\in L^2(0,T;L^2(\Omega_p)),\\
p\in L^2(0,T;H^{j+1}(\Omega_f)).
\end{cases}
\end{equation}
And we suppose that
\begin{align}\label{regularity right}
u_h^0\in L^2(\Omega_f),\quad
\phi_h^0\in L^2(\Omega_p),\quad
f_1\in L^2(0,T;L^2(\Omega_f)),
\quad f_2\in L^2(0,T;L^2(\Omega_p)).
\end{align}

\subsection{First-order scheme}
By using the backward Euler method, we obtain the first-order scheme in time.\\
\textbf{Algorithm 1} (First-order Scheme)
Given $(u_{h}^n,p_{h}^n,\phi_{h}^n, r_h^n )$, for any $(v_h,q_h,\psi_h)\in H_{fh}\times Q_{fh}\times H_{ph}$, find $(u_{h}^{n+1},p_{h}^{n+1},\phi_{h}^{n+1}, r_h^{n+1} )$ satisfying
\begin{align}
&(\frac{u_{h}^{n+1}-u_{h}^n}{\Delta t},v_h)_f+\nu(\nabla u_{h}^{n+1},\nabla v_h)_f+{\eta}\sum_{i=1}^{d-1}\int_{\Gamma}(u_{h}^{n+1}\cdot{\tau}_i)\cdot(v_h\cdot {\tau}_i)ds
-(p_{h}^{n+1},\nabla\cdot v_h)_f\nonumber\\
&+\frac{r_{h}^{n+1}}{exp(-\frac{t^{n+1}}{T})}a_N(u_h^n,u_h^n,v_h)
+\frac{r_{h}^{n+1}}{exp(-\frac{t^{n+1}}{T})}c_\Gamma(v_h,\phi_{h}^n)
=(f_1^{n+1},v_h)_f,\label{A1stokes_h1}\\
&(\nabla\cdot u_{h}^{n+1},q_h)_f=0,\label{A1stokes_h2}\\
&gS_0(\frac{\phi_{h}^{n+1}-\phi_{h}^n}{\Delta t},\psi_h)_p
+g({\mathcal{K}}\nabla\phi_{h}^{n+1},\nabla\psi_h)_p
-\frac{r_{h}^{n+1}}{exp(-\frac{t^{n+1}}{T})}c_\Gamma(u_{h}^n,\psi_h)
=g(f_2^{n+1},\psi_h)_p,\label{A1darcy_h}\\
&\frac{r_{h}^{n+1}-r_{h}^n}{\Delta t}=-\frac{1}{T}r_{h}^{n+1}+\frac{1}{exp(-\frac{t^{n+1}}{T})}
(c_\Gamma(u_{h}^{n+1},\phi_{h}^n)-c_\Gamma(u_{h}^n,\phi_{h}^{n+1})
+a_N(u_h^n,u_h^n,u_h^{n+1})).\label{A1r}
\end{align}
\medskip
\begin{theorem}\label{stability of algo1}
For any $N\geq1$, Algorithm 1 satisfies the following unconditional stability properties
\begin{equation}\label{A1 stable}
\begin{split}
&\|u_h^N\|_{0,f}^2
+\sum_{n=0}^{N-1}\|u_h^{n+1}-u_h^n\|_{0,f}^2
+\nu\Delta t\sum_{n=0}^{N-1}\| \nabla u_h^{n+1}\|_{0,f}^2
+\eta\Delta t\sum_{i=1}^{d-1}\int_\Gamma(u_h^N\cdot \tau_i)^2ds
+{gS_0}\|\phi_h^N\|_{0,p}^2\nonumber\\
&
+gS_0\sum_{n=0}^{N-1}\|\phi_h^{n+1}-\phi_h^n\|_{0,p}^2
+gk\Delta t\sum_{n=0}^{N-1}\|\nabla \phi_h^{n+1}\|_{0,p}^2
+|r_h^N|^2
+\sum_{n=0}^{N-1}|r_h^{n+1}-r_h^n|^2
+\frac{2\Delta t}{T}\sum_{n=0}^{N-1}|r_h^{n+1}|^2\\
&\leq C\Delta t\sum_{n=0}^{N-1}\|f_1^{n+1}\|_{0,f}^2
+C\Delta t\sum_{n=0}^{N-1}\|f_2^{n+1}\|_{0,p}^2
+\|u_h^0\|_{0,f}^2
+\eta\Delta t\sum_{i=1}^{d-1}\int_\Gamma (u_h^0\cdot{\tau}_i)^2ds
+gS_0\|\phi_h^0\|_{0,p}^2
+|r_h^0|^2,
\end{split}
\end{equation}
where the positive constant $C$ is independent of $h$ and $\Delta t$.
\end{theorem}
%
\medskip
\begin{proof}
Substituting $v_h=u_h^{n+1},~\psi_h=\phi_h^{n+1}$ into (\ref{A1stokes_h1}),(\ref{A1darcy_h}) and combining with (\ref{A1stokes_h2}), yields
\begin{align}
&(\frac{u_h^{n+1}-u_h^n}{\Delta t},u_h^{n+1})_f
+\nu(\nabla u_h^{n+1}, \nabla u_h^{n+1})_{f}
+\eta\sum_{i=1}^{d-1}\int_\Gamma(u_{h}^{n+1}\cdot{\tau}_i)^2ds
+\frac{r_h^{n+1}}{exp(-\frac{t^{n+1}}{T})}a_N(u_h^n,u_h^n,u_h^{n+1})\nonumber\\
&+\frac{r_h^{n+1}}{exp(-\frac{t^{n+1}}{T})}c_\Gamma(u_h^{n+1},\phi_h^n)
=(f_1^{n+1},u_h^{n+1})_f,\label{A1s_uu}\\
&gS_0(\frac{\phi_h^{n+1}-\phi_h^n}{\Delta t},\phi_h^{n+1})_p
+g(\mathcal{K}\nabla\phi_h^{n+1},\nabla\phi_h^{n+1})_p
-\frac{r_h^{n+1}}{exp(-\frac{t^{n+1}}{T})}c_\Gamma(u_h^n,\phi_h^{n+1})=g(f_2^{n+1},\phi_h^{n+1})_p.\label{A1d_phiphi}
\end{align}
Multiplying (\ref{A1r}) by $r_h^{n+1}$, we have
\begin{align}\label{A1rr}
\frac{r_{h}^{n+1}-r_{h}^n}{\Delta t} r_h^{n+1} =-\frac{1}{T}|r_{h}^{n+1}|^2+\frac{r_h^{n+1}}{exp(-\frac{t^{n+1}}{T})}
(c_\Gamma(u_{h}^{n+1},\phi_{h}^n)-c_\Gamma(u_{h}^n,\phi_{h}^{n+1})
+a_N(u_h^n,u_h^n,u_h^{n+1})).
\end{align}
Adding (\ref{A1s_uu}), (\ref{A1d_phiphi}) and (\ref{A1rr}) together and using the fact that
$(a-b,a)=\frac{1}{2}(|a|^2-|b|^2+|a-b|^2),$ we obtain
\begin{align}\label{A1after-a-b a}
&\frac{1}{2\Delta t}\|u_h^{n+1}\|_{0,f}^2
-\frac{1}{2\Delta t}\|u_h^n\|_{0,f}^2
+\frac{1}{2\Delta t}\|u_h^{n+1}-u_h^n\|_{0,f}^2
+\nu\|\nabla u_h^{n+1}\|_{0,f}^2
+\eta\sum_{i=1}^{d-1}\int_\Gamma(u_{h}^{n+1}\cdot{\tau}_i)^2ds\nonumber\\
&
+\frac{gS_0}{2\Delta t}\|\phi_h^{n+1}\|_{0,p}^2
-\frac{gS_0}{2\Delta t}\|\phi_h^n\|_{0,p}^2
+\frac{gS_0}{2\Delta t}\|\phi_h^{n+1}-\phi_h^n\|_{0,p}^2
+gk\|\nabla \phi_h^{n+1}\|_{0,p}^2
+\frac{1}{2\Delta t}|r_h^{n+1}|^2\nonumber\\
&
-\frac{1}{2\Delta t}|r_h^n|^2+\frac{1}{2\Delta t}|r_h^{n+1}-r_h^n|^2+\frac{1}{T}|r_h^{n+1}|^2
=(f_1^{n+1},u_h^{n+1})_f+g(f_2^{n+1},\phi_h^{n+1})_p.
\end{align}
Next we can estimate the right hand side terms of \eqref{A1after-a-b a} by using Young's and Poincar${\rm \acute{e}}$ inequalities:
\begin{align}\label{A1right}
&(f_1^{n+1},u_h^{n+1})_f+g(f_2^{n+1},\phi_h^{n+1})_p\nonumber\\
&\leq C\|f_1^{n+1}\|_{0,f}\|u_h^{n+1}\|_{0,f}
+C\|f_2^{n+1}\|_{0,p}\|\phi_h^{n+1}\|_{0,p}\nonumber\\
&\leq C\|f_1^{n+1}\|_{0,f}^2+\frac{\nu}{2}\|\nabla u_h^{n+1}\|_{0,f}^2
+C\|f_2^{n+1}\|_{0,p}^2+\frac{gk}{2}\|\nabla \phi_h^{n+1}\|_{0,p}^2.
\end{align}
Combining (\ref{A1after-a-b a}) and (\ref{A1right}) gives that
\begin{align}\label{A1u n+1 - u n}
&\frac{1}{2\Delta t}\|u_h^{n+1}\|_{0,f}^2-\frac{1}{2\Delta t}\|u_h^n\|_{0,f}^2
+\frac{1}{2\Delta t}\|u_h^{n+1}-u_h^n\|_{0,f}^2
+\frac{\nu}{2}\|\nabla u_h^{n+1}\|_{0,f}^2
+\frac{\eta}{2}\sum_{i=1}^{d-1}\int_\Gamma (u_h^{n+1}\cdot \tau_i)^2ds\nonumber\\
&-\frac{\eta}{2}\sum_{i=1}^{d-1}\int_\Gamma (u_h^n\cdot \tau_i)^2ds
+\frac{gS_0}{2\Delta t}\|\phi_h^{n+1}\|_{0,p}^2
-\frac{gS_0}{2\Delta t}\|\phi_h^n\|_{0,p}^2
+\frac{gS_0}{2\Delta t}\|\phi_h^{n+1}-\phi_h^n\|_{0,p}^2
+\frac{gk}{2}\|\nabla \phi_h^{n+1}\|_{0,p}^2\nonumber\\
&
+\frac{1}{2\Delta t}|r_h^{n+1}|^2-\frac{1}{2\Delta t}|r_h^n|^2+\frac{1}{2\Delta t}|r_h^{n+1}-r_h^n|^2
+\frac{1}{T}|r_h^{n+1}|^2
\leq C\|f_1^{n+1}\|_{0,f}^2+C\|f_2^{n+1}\|_{0,p}^2.
\end{align}
Summing up (\ref{A1u n+1 - u n}) from $n=0$ to $n=N-1$, and multiplying by $2\Delta t$, we get
\begin{align*}
&\|u_h^N\|_{0,f}^2
+\sum_{n=0}^{N-1}\|u_h^{n+1}-u_h^n\|_{0,f}^2
+\nu\Delta t\sum_{n=0}^{N-1}\|\nabla u_h^{n+1}\|_{0,f}^2
+\eta\Delta t\sum_{i=1}^{d-1}\int_\Gamma(u_h^N\cdot \tau_i)^2ds
+{gS_0}\|\phi_h^N\|_{0,p}^2\nonumber\\
&
+gS_0\sum_{n=0}^{N-1}\|\phi_h^{n+1}-\phi_h^n\|_{0,p}^2
+gk\Delta t\sum_{n=0}^{N-1}\|\nabla \phi_h^{n+1}\|_{0,p}^2
+|r_h^N|^2
+\sum_{n=0}^{N-1}|r_h^{n+1}-r_h^n|^2
+\frac{2\Delta t}{T}\sum_{n=0}^{N-1}|r_h^{n+1}|^2\\
&\leq C\Delta t\sum_{n=0}^{N-1}\|f_1^{n+1}\|_{0,f}^2
+C\Delta t\sum_{n=0}^{N-1}\|f_2^{n+1}\|_{0,p}^2
+\|u_h^0\|_{0,f}^2
+\eta\Delta t\sum_{i=1}^{d-1}\int_\Gamma (u_h^0\cdot{\tau}_i)^2ds
+gS_0\|\phi_h^0\|_{0,p}^2
+|r_h^0|^2.
\end{align*}
\end{proof}\\

\noindent\textbf{Remark 1.} By using Theorem \ref{stability of algo1}, we can easily obtain that
\begin{align}\label{uhn1}
\Delta t \sum_{n=0}^{N-1} \|\nabla u_h^{n+1}\|_{0,f}^2\leq C^*,
\end{align}
where the positive constant $C^*$ is independent of $h$ and $\Delta t$. The boundedness of the velocity is particularly useful in the error analysis.

\subsection{Second-order scheme}
By using the second-order backward differentiation formulation (BDF2), we obtain the second-order scheme in time.\\
\textbf{Algorithm 2} (Second-order Scheme)
Given $(u_{h}^{n-1},p_{h}^{n-1},\phi_{h}^{n-1}, r_h^{n-1} )$ and $(u_{h}^n,p_{h}^n,\phi_{h}^n, r_h^n)$, then for any $(v_h,q_h,\psi_h)\in H_{fh}\times Q_{fh}\times H_{ph}$, find $(u_{h}^{n+1},p_{h}^{n+1},\phi_{h}^{n+1}, r_h^{n+1})$ satisfying
\begin{align}
&(\frac{3u_h^{n+1}-4u_h^n+u_h^{n-1}}{2\Delta t},v_h)_f
+\nu(\nabla u_h^{n+1},\nabla v_h)_f
+\eta\sum_{i=1}^{d-1}\int_\Gamma(u_h^{n+1}\cdot{\tau}_i)(v_h\cdot{\tau}_i)ds
-(p_h^{n+1},\nabla\cdot v_h)_f\nonumber\\
&
+\frac{r_h^{n+1}}{exp(-\frac{t^{n+1}}{T})}a_N(2u_h^n-u_h^{n-1},2u_h^n-u_h^{n-1},v_h)
+\frac{r_h^{n+1}}{exp(-\frac{t^{n+1}}{T})}c_\Gamma(v_h,2\phi_h^n-\phi_h^{n-1})
=(f_1^{n+1},v_h)_f,\label{A2stokes1}\\
&(\nabla\cdot u_h^{n+1},q_h)_f=0,\label{A2stokes2}\\
&gS_0(\frac{3\phi_h^{n+1}-4\phi_h^n+\phi_h^{n-1}}{2\Delta t},\psi_h)_p
+g(\mathcal{K}\nabla\phi_h^{n+1},\nabla\psi_h)_p
-\frac{r_h^{n+1}}{exp(-\frac{t^{n+1}}{T})}c_\Gamma(2u_h^n-u_h^{n-1},\psi_h)
=g(f_2^{n+1},\psi_h)_p,\label{A2darcy}\\
&\frac{3r_h^{n+1}-4r_h^n+r_h^{n-1}}{2\Delta t}
=-\frac{1}{T}r_h^{n+1}
+\frac{1}{exp(-\frac{t^{n+1}}{T})}(c_\Gamma(u_h^{n+1},2\phi_h^n-\phi_h^{n-1})
-c_\Gamma(2u_h^n-u_h^{n-1},\phi_h^{n+1})\nonumber\\
&+a_N(2u_h^n-u_h^{n-1},2u_h^n-u_h^{n-1},u_h^{n+1})).\label{A2r}
\end{align}

\medskip
\begin{theorem}
For any $N\geq2$, Algorithm 2 satisfies the following unconditional stability properties
\begin{equation}\label{Algoritgm 2 stability}\begin{split}
&\|u_h^{N}\|_{0,f}^2
+\|2u_h^{N}-u_h^{N-1}\|_{0,f}^2
+\sum_{n=1}^{N-1}\|u_h^{n+1}-2u_h^n+u_h^{n-1}\|_{0,f}^2
+2\nu\Delta t\sum_{n=1}^{N-1}\|\nabla u_h^{n+1}\|_{0,f}\\
&
+{2\eta}\Delta t\sum_{i=1}^{d-1}\int_\Gamma(u_h^N\cdot{\tau}_i)^2ds
+gS_0\|\phi_h^N\|_{0,p}^2
+gS_0\|2\phi_h^N-\phi_h^{N-1}\|_{0,p}^2
+{gS_0}\sum_{n=1}^{N-1}\|\phi_h^{n+1}-2\phi_h^n+\phi_h^{n-1}\|_{0,p}^2\\
&
+2gk\Delta t\sum_{n=1}^{N-1}\|\nabla \phi_h^{n+1}\|_{0,p}^2
+|r_h^{N}|^2
+|2r_h^{N}-r_h^{N-1}|^2
+\sum_{n=1}^{N-1}|r_h^{n+1}-2r_h^n+r_h^{n-1}|^2
+\frac{4\Delta t}{T}\sum_{n=1}^{N-1}|r_h^{n+1}|^2\\
&\leq C\Delta t\sum_{n=1}^{N-1}\|f_1^{n+1}\|_{0,f}^2
+C \Delta t\sum_{n=1}^{N-1}\|f_2^{n+1}\|_{0,p}^2
+\|u_h^{1}\|_{0,f}^2
+\|2u_h^{1}-u_h^0\|_{0,f}^2
+{2\eta}\Delta t\sum_{i=1}^{d-1}\int_\Gamma(u_h^1\cdot{\tau}_i)^2ds\\
&
+{gS_0}\|\phi_h^1\|_{0,p}^2
+{gS_0}\|2\phi_h^1-\phi_h^{0}\|_{0,p}^2
+|r_h^{1}|^2
+|2r_h^{1}-r_h^{0}|^2,
\end{split}
\end{equation}
where the positive constant $C$ is independent of $h$ and $\Delta t$.
\end{theorem}

\medskip
\begin{proof}
Multiplying (\ref{A2r}) by $r_h^{n+1}$ yields
\begin{align}\label{A2rr}
&\frac{3r_h^{n+1}-4r_h^n+r_h^{n-1}}{2\Delta t} r_h^{n+1}
=-\frac{1}{T}|r_h^{n+1}|^2
+\frac{r_h^{n+1}}{exp(-\frac{t^{n+1}}{T})}(c_\Gamma(u_h^{n+1},2\phi_h^n-\phi_h^{n-1})
-c_\Gamma(2u_h^n-u_h^{n-1},\phi_h^{n+1})\nonumber\\
&+a_N(2u_h^n-u_h^{n-1},2u_h^n-u_h^{n-1},u_h^{n+1})).
\end{align}
Setting $v_h=u_h^{n+1},~\psi_h=\phi_h^{n+1}$ in (\ref{A2stokes1})-(\ref{A2darcy}), and adding these equations with (\ref{A2rr}) together, we have
\begin{align}\label{A2 341}
&(\frac{3u_h^{n+1}-4u_h^n+u_h^{n-1}}{2\Delta t},u_h^{n+1})_f
+\nu\|\nabla u_h^{n+1}\|_{0,f}^2
+\eta\sum_{i=1}^{d-1}\int_\Gamma(u_h^{n+1}\cdot{\tau}_i)^2ds\nonumber\\
&+gS_0(\frac{3\phi_h^{n+1}-4\phi_h^n+\phi_h^{n-1}}{2\Delta t},\phi_h^{n+1})_p
+gk\|\nabla \phi_h^{n+1}\|_{0,p}^2
+\frac{3r_h^{n+1}-4r_h^n+r_h^{n-1}}{2\Delta t} r_h^{n+1}
+\frac{1}{T}|r_h^{n+1}|^2\nonumber\\
&
=(f_1^{n+1},u_h^{n+1})_f+g(f_2^{n+1},\phi_h^{n+1})_p.
\end{align}
Using the fact that $(3a-4b+c,a)=\frac{1}{2}(|a|^2+|2a-b|^2-|b|^2-|2b-c|^2+|a-2b+c|^2)$, we can finally obtain (\ref{Algoritgm 2 stability}) with the similar procedure as in \cite{jiang2021SAV}.
\end{proof}

\section{Error Analysis}
In this section, we will conduct a rigorous analysis of the error convergence rates of the first-order unconditionally stable SAV method for the unsteady Navier-Stokes-Darcy problem.

We introduce the Ritz projections $\Pi_h:~H_f\rightarrow H_{fh}$ such that
\begin{align*}
(\nabla(v-\Pi_h v),\nabla v_h)_f=0,\quad\forall v\in H_f,~v_h\in H_{fh},
\end{align*}
and $R_h:~H_p\rightarrow H_{ph}$ such that
\begin{align*}
(\nabla(\psi-R_h\psi),\nabla \psi_h)_p=0,\quad\forall \psi\in H_p,~\psi_h\in H_{ph},
\end{align*}
then for $2\leq k\leq 6$, there holds
\begin{equation}\label{Ritz u}
\begin{cases}
  \|v-\Pi_h v\|_{L^k(\Omega_f)}+h\|\nabla(v-\Pi_h v)\|_{L^k(\Omega_f)}
 \leq  Ch^{r+1}\|v\|_{W^{r+1,k}(\Omega_f)},\quad\forall v\in H_f,\\
 \|\frac{\partial}{\partial t}(v^n-\Pi_h v^n)\|_{L^2(\Omega_f)}
 +h\|\nabla(\frac{\partial}{\partial t}(v^n-\Pi_h v^n))\|_{L^2(\Omega_f)}
\leq Ch^{r+1}\|\frac{\partial v^n}{\partial t}\|_{H^{r+1}(\Omega_f)}
,
\end{cases}
\end{equation}
and
\begin{equation}\label{Ritz phi}
\begin{cases}
  \|\psi-R_h \psi\|_{L^k(\Omega_p)}+h\|\nabla(\psi-R_h \psi)\|_{L^k(\Omega_p)}
 \leq  Ch^{r+1}\|\psi\|_{W^{r+1,k}(\Omega_p)},\quad\forall \psi\in H_p,\\
 \|\frac{\partial}{\partial t}(\psi^n-R_h \psi^n)\|_{L^2(\Omega_p)}
 +h\|\nabla(\frac{\partial}{\partial t}(\psi^n-R_h \psi^n))\|_{L^2(\Omega_p)}
\leq Ch^{r+1}\|\frac{\partial \psi^n}{\partial t}\|_{H^{r+1}(\Omega_p)}
.
\end{cases}
\end{equation}

We first split the error as:
\begin{align*}
&e_u^{n+1}=u^{n+1}-u_h^{n+1}=u^{n+1}-\Pi_h u^{n+1}+\Pi_h u^{n+1}-u_h^{n+1}=\sigma_u^{n+1}+\theta_u^{n+1},\\
&e_\phi^{n+1}=\phi^{n+1}-\phi_h^{n+1}=\phi^{n+1}-R_h \phi^{n+1}+R_h \phi^{n+1}-\phi_h^{n+1}
=\sigma_\phi^{n+1}+\theta_\phi^{n+1},\\
&e_r^{n+1}=r^{n+1}-r_h^{n+1}.
\end{align*}

\begin{theorem} \label{final thm}
Under the regularity assumptions (\ref{regularity solution}) and (\ref{regularity right}), the numerical solution obtained from Algorithm 1 satisfies
\begin{align}\label{err 1order}
&\|e_u^s\|_{0,f}^2
+\sum_{n=0}^{s-1}\|e_u^{n+1}-e_u^n\|_{0,f}^2
+\nu\Delta t\sum_{n=0}^{s-1}\|\nabla e_u^{n+1}\|_{0,f}^2
+\frac{4\nu}{9}\Delta t\|\nabla e_u^s\|_{0,f}^2
+\frac{\eta}{2}\Delta t\sum_{i=1}^{d-1}\int_\Gamma(e_u^s\cdot\tau_i)^2ds\nonumber\\
&+gS_0\|e_\phi^s\|_{0,p}^2
+gS_0\sum_{n=0}^{s-1}\|e_\phi^{n+1}-e_\phi^n\|_{0,p}^2
+gk\Delta t\sum_{n=0}^{s-1}\|\nabla e_\phi^{n+1}\|_{0,p}^2
+\frac{2gk}{5}\Delta t\|\nabla e_\phi^s\|_{0,p}^2
+\frac{1}{2}|e_r^s|^2\nonumber\\
&+\sum_{n=0}^{s-1}|e_r^{n+1}-e_r^n|^2
+\frac{\Delta t}{T}\sum_{n=0}^{s-1}|e_r^{n+1}|^2\nonumber\\
&\leq
C\Delta t^2
(1+\|\partial_t u\|^2_{2,0,f}+\|\partial_t u\|^2_{2,1,f}+\|\partial_t u\|^2_{2,2,f}
+\|\partial_{tt} u\|^2_{2,0,f}+\|\partial_t \phi\|^2_{2,0,p}
+\|\partial_t \phi\|^2_{2,1,p}+\|\partial_{tt} \phi\|^2_{2,0,p})\nonumber\\
&\quad
+Ch^{2k}\normmm{u}^2_{2,k+1,f}
+Ch^{2k}\normmm{u}^2_{\infty,k+1,f}
+Ch^{2k+2}\|\partial_t u\|^2_{2,k+1,f}
+Ch^{2m}\normmm{\phi}^2_{2,m+1,p}\nonumber\\
&\quad
+Ch^{2m+2}\|\partial_t \phi\|^2_{2,m+1,p}
+Ch^{2k}\normmm{u}_{2,k+1,f}^2
+Ch^{2m+2}\normmm{\phi}_{\infty,m+1,p}^2
+Ch^{2m}\normmm{\phi}_{2,m+1,p}^2,
\end{align}
where the positive constant $C$ is independent of $h$ and $\Delta t$, and $1\leq s\leq N$.
\end{theorem}

\medskip
\begin{proof} We divide the error estimate into the following four steps to obtain the desired result.

\textbf{ Step 1. Determination of error equations}

For any $v_h\in V_h,~\psi_h\in H_{ph}$, the exact solution $(u^{n+1},p^{n+1},\phi^{n+1})$ of problem (\ref{nsd problem}) satisfies
\begin{align}
&(\frac{u^{n+1}-u^n}{\Delta t},v_h)_f
+\nu(\nabla u^{n+1},\nabla v_h)_f
+\eta\sum_{i=1}^{d-1}\int_\Gamma(u^{n+1}\cdot\tau_i)(v_h\cdot\tau_i)ds
-(p^{n+1},\nabla\cdot v_h)_f\label{exact1}\\
&+a_N(u^{n+1},u^{n+1},v_h)
+c_\Gamma(v_h,\phi^{n+1})
=(f_1^{n+1},v_h)_f+R_u^{n+1}(v_h),\nonumber\\
&gS_0(\frac{\phi^{n+1}-\phi^n}{\Delta t},\psi_h)_p
+g(\mathcal{K}\nabla\phi^{n+1},\nabla\psi_h)_p
-c_\Gamma(u^{n+1},\psi_h)
=g(f_2^{n+1},\psi_h)_p+R_\phi^{n+1}(\psi_h).\label{exact2}
\end{align}
And we have
\begin{align}
&\frac{r^{n+1}-r^n}{\Delta t}=-\frac{1}{T}r^{n+1}+\frac{1}{exp(-\frac{t^{n+1}}{T})}
(a_N(u^{n+1},u^{n+1},u^{n+1})+c_\Gamma(u^{n+1},\phi^{n+1})\nonumber\\
&\qquad\qquad\qquad-c_\Gamma(u^{n+1},\phi^{n+1}))
+R_r^{n+1}.\label{exact3}
\end{align}
The definitions of above error terms are given by
\begin{align*}
  R_u^{n+1}(v_h)&=(\frac{u^{n+1}-u^n}{\Delta t}-\partial _t u(t^{n+1}),v_h)_f, \\
  R_\phi^{n+1}(\psi_h)&=gS_0(\frac{\phi^{n+1}-\phi^n}{\Delta t}-\partial_t \phi(t^{n+1}),\psi_h)_p, \\
  R_r^{n+1}&=\frac{r^{n+1}-r^n}{\Delta t}-r_t^{n+1}.
\end{align*}
Subtracting (\ref{A1stokes_h1})-(\ref{A1r}) from (\ref{exact1})-(\ref{exact3}), we obtain
\begin{align}
&(\frac{e_u^{n+1}-e_u^n}{\Delta t},v_h)_f
+\nu(\nabla e_u^{n+1},\nabla v_h)_f
+\eta\sum_{i=1}^{d-1}
\int_\Gamma(e_u^{n+1}\cdot \tau_i)(v_h\cdot\tau_i)ds
-(p^{n+1}-p_h^{n+1},\nabla\cdot v_h)_f\nonumber\\
&+a_N(u^{n+1},u^{n+1},v_h)-\frac{r_h^{n+1}}{exp \ttt) }a_N(u_h^n,u_h^n,v_h)
+c_\Gamma(v_h,\phi^{n+1})\nonumber\\
&
-\frac{r_h^{n+1}}{exp\ttt)}c_\Gamma(v_h,\phi_h^n)=R_u^{n+1}(v_h),\label{err1}\\
&gS_0(\frac{e_\phi^{n+1}-e_\phi^n}{\Delta t},\psi_h)_p
+g(\mathcal{K}\nabla e_\phi^{n+1},\nabla \psi_h)_p
-c_\Gamma(u^{n+1},\psi_h)+\frac{r_h^{n+1}}{exp\ttt)}c_\Gamma(u_h^n,\psi_h)
=R_\phi^{n+1}(\psi_h),\label{err2}\\
&\frac{e_r^{n+1}-e_r^n}{\Delta t}+\frac{1}{T}e_r^{n+1}=\frac{1}{exp\ttt)}(a_N(u^{n+1},u^{n+1},u^{n+1})
-a_N(u_h^n,u_h^n,u_h^{n+1})+c_\Gamma(u^{n+1},\phi^{n+1})\nonumber\\
&-c_\Gamma(u_h^{n+1},\phi_h^n)
-c_\Gamma(u^{n+1},\phi^{n+1})+c_\Gamma(u_h^n,\phi_h^{n+1}))+R_r^{n+1}.\label{err3}
\end{align}
Letting $v_h=\theta _u^{n+1}$, $\psi_h=\theta_\phi^{n+1}$ in (\ref{err1}) and (\ref{err2}), we have
\begin{align}
&(\frac{\theta_u^{n+1}-\theta_u^n}{\Delta t},\theta_u^{n+1})_f
+\nu(\nabla \theta_u^{n+1},\nabla \theta_u^{n+1})_f
+\eta\sum_{i=1}^{d-1}\int_\Gamma(\theta_u^{n+1}\cdot\tau_i)^2ds
+a_N(u^{n+1},u^{n+1},\theta_u^{n+1})\nonumber\\
&
-\frac{r_h^{n+1}}{exp\ttt)}a_N(u_h^n,u_h^n,\theta_u^{n+1})
+c_\Gamma(\theta_u^{n+1},\phi^{n+1})
-\frac{r_h^{n+1}}{exp\ttt)}c_\Gamma(\theta_u^{n+1},\phi_h^n)\nonumber\\
&\quad=
R_u^{n+1}(\theta_u^{n+1})-(\frac{\sigma_u^{n+1}-\sigma_u^n}{\Delta t},\theta_u^{n+1})_f
-\eta\sum_{i=1}^{d-1}\int_\Gamma(\sigma_u^{n+1}\cdot\tau_i)
(\theta_u^{n+1}\cdot\tau_i)ds,\label{ns1}\\
&gS_0(\frac{\theta_\phi^{n+1}-\theta_\phi^n}{\Delta t},\theta_\phi^{n+1})_p
+g(\mathcal{K}\nabla\theta_\phi^{n+1},\nabla\theta_\phi^{n+1})_p
-c_\Gamma(u^{n+1},\theta_\phi^{n+1})
+\frac{r_h^{n+1}}{exp\ttt)}c_\Gamma(u_h^n,\theta_\phi^{n+1})\nonumber\\
&\quad=R_\phi^{n+1}(\theta_\phi^{n+1})
-gS_0(\frac{\sigma_\phi^{n+1}-\sigma_\phi^n}{\Delta t},\theta_\phi^{n+1})_p
.\label{d2}
\end{align}
Multiplying (\ref{err3}) by $e_r^{n+1}$ on both sides, we derive
\begin{align}
&\frac{e_r^{n+1}-e_r^n}{\Delta t} e_r^{n+1}+\frac{1}{T}|e_r^{n+1}|^2=
\frac{e_r^{n+1}}{exp\ttt)}(a_N(u^{n+1},u^{n+1},u^{n+1})-a_N(u_h^n,u_h^n,u_h^{n+1})
\nonumber\\
&\qquad\qquad\qquad\qquad\qquad\qquad~~
+c_\Gamma(u^{n+1},\phi^{n+1})
-c_\Gamma(u_h^{n+1},\phi_h^n)-c_\Gamma(u^{n+1},\phi^{n+1})
\nonumber\\
&\qquad\qquad\qquad\qquad\qquad\qquad~~
+c_\Gamma(u_h^n,\phi_h^{n+1}))
+R_r^{n+1} e_r^{n+1}.\label{er3}
\end{align}
Using the fact that $(a-b,a)=\frac{1}{2}(|a|^2-|b|^2+|a-b|^2),$ and adding (\ref{ns1}) and (\ref{d2}) together, we have
\begin{align}\label{err theta}
&\frac{1}{2\Delta t}\|\theta_u^{n+1}\|_{0,f}^2
-\frac{1}{2\Delta t}\|\theta_u^{n}\|_{0,f}^2
+\frac{1}{2\Delta t}\|\theta_u^{n+1}-\theta_u^n\|_{0,f}^2
+\nu\|\nabla\theta_u^{n+1}\|_{0,f}^2
+\eta\sum_{i=1}^{d-1}\int_\Gamma(\theta_u^{n+1}\cdot\tau_i)^2ds\nonumber\\
&+\frac{gS_0}{2\Delta t}\|\theta_\phi^{n+1}\|_{0,p}^2
-\frac{gS_0}{2\Delta t}\|\theta_\phi^{n}\|_{0,p}^2
+\frac{gS_0}{2\Delta t}\|\theta_\phi^{n+1}-\theta_\phi^n\|_{0,p}^2
+gk\|\nabla\theta_\phi^{n+1}\|_{0,p}^2\nonumber\\
&
=
-a_N(u^{n+1},u^{n+1},\theta_u^{n+1})
+\frac{r_h^{n+1}}{exp\ttt)}a_N(u_h^n,u_h^n,\theta_u^{n+1})
-c_\Gamma(\theta_u^{n+1},\phi^{n+1})
\nonumber\\
&\quad
+\frac{r_h^{n+1}}{exp\ttt)}c_\Gamma(\theta_u^{n+1},\phi_h^n)
+c_\Gamma(u^{n+1},\theta_\phi^{n+1})
-\frac{r_h^{n+1}}{exp\ttt)}c_\Gamma(u_h^n,\theta_\phi^{n+1})
\nonumber\\
&\quad+R_u^{n+1}(\theta_u^{n+1})+R_\phi^{n+1}(\theta_\phi^{n+1})
-(\frac{\sigma_u^{n+1}-\sigma_u^n}{\Delta t},\theta_u^{n+1})_f-gS_0(\frac{\sigma_\phi^{n+1}-\sigma_\phi^n}{\Delta t},\theta_\phi^{n+1})_p
\nonumber\\
&\quad
-\eta\sum_{i=1}^{d-1}\int_\Gamma(\sigma_u^{n+1}\cdot\tau_i)
(\theta_u^{n+1}\cdot\tau_i)ds.
\end{align}
We rewrite the trilinear terms and coupling terms on the right hand side of \eqref{err theta} as
\begin{flalign*}
&{\bm {T_1}}=-a_N(u^{n+1},u^{n+1},\theta_u^{n+1})+\frac{r_h^{n+1}}{exp\ttt)}a_N(u_h^n,u_h^n,\theta_u^{n+1})\\
&\quad=-a_N(u^{n+1}-u_h^n,u^{n+1},\theta_u^{n+1})-a_N(u_h^n,u^{n+1},\theta_u^{n+1})
+a_N(u_h^n,u_h^n,\theta_u^{n+1})-a_N(u_h^n,u_h^n,\theta_u^{n+1})\\
&\quad\quad+\frac{r_h^{n+1}}{exp\ttt)}a_N(u_h^n,u_h^n,\theta_u^{n+1})\\
&\quad=-a_N(u^{n+1}-u_h^n,u^{n+1},\theta_u^{n+1})-a_N(u_h^n,u^{n+1}-u_h^n,\theta_u^{n+1})
-\frac{e_r^{n+1}}{exp\ttt)}a_N(u_h^n,u_h^n,\theta_u^{n+1}),&
\end{flalign*}

\begin{flalign*}
&{\bm{T_2}}=-c_\Gamma(\theta_u^{n+1},\phi^{n+1})
+\frac{r_h^{n+1}}{exp\ttt)}c_\Gamma(\theta_u^{n+1},\phi_h^n)\\
&\quad=-c_\Gamma(\theta_u^{n+1},\phi^{n+1}-\phi_h^n)-c_\Gamma(\theta_u^{n+1},\phi_h^n)
+\frac{r_h^{n+1}}{exp\ttt)}c_\Gamma(\theta_u^{n+1},\phi_h^n)\\
&\quad=-c_\Gamma(\theta_u^{n+1},\phi^{n+1}-\phi_h^n)
-\frac{e_r^{n+1}}{exp\ttt)}c_\Gamma(\theta_u^{n+1},\phi_h^n),&
\end{flalign*}
\begin{flalign*}
&{\bm{T_3}}=c_\Gamma(u^{n+1},\theta_\phi^{n+1})
-\frac{r_h^{n+1}}{exp\ttt)}c_\Gamma(u_h^n,\theta_\phi^{n+1})\\
&\quad=c_\Gamma(u^{n+1}-u_h^n,\theta_\phi^{n+1})+c_\Gamma(u_h^n,\theta_\phi^{n+1})
-\frac{r_h^{n+1}}{exp\ttt)}c_\Gamma(u_h^n,\theta_\phi^{n+1})\\
&\quad=c_\Gamma(u^{n+1}-u_h^n,\theta_\phi^{n+1})
+\frac{e_r^{n+1}}{exp\ttt)}c_\Gamma(u_h^n,\theta_\phi^{n+1}).&
\end{flalign*}
%
%
%
Hence (\ref{err theta}) can be recast as
\begin{flalign}
&\frac{1}{2\Delta t}\|\theta_u^{n+1}\|_{0,f}^2
-\frac{1}{2\Delta t}\|\theta_u^{n}\|_{0,f}^2
+\frac{1}{2\Delta t}\|\theta_u^{n+1}-\theta_u^n\|_{0,f}^2
+\nu\|\nabla\theta_u^{n+1}\|_{0,f}^2
+\eta\sum_{i=1}^{d-1}\int_\Gamma(\theta_u^{n+1}\cdot\tau_i)^2ds \nonumber \\
&+\frac{gS_0}{2\Delta t}\|\theta_\phi^{n+1}\|_{0,p}^2
-\frac{gS_0}{2\Delta t}\|\theta_\phi^{n}\|_{0,p}^2
+\frac{gS_0}{2\Delta t}\|\theta_\phi^{n+1}-\theta_\phi^n\|_{0,p}^2
+gk\|\nabla\theta_\phi^{n+1}\|_{0,p}^2 \label{err ADD}  \\
&=
-a_N(u^{n+1}-u_h^n,u^{n+1},\theta_u^{n+1})
-a_N(u_h^n,u^{n+1}-u_h^n,\theta_u^{n+1})\nonumber\\
&\quad
-c_\Gamma(\theta_u^{n+1},\phi^{n+1}-\phi_h^n)
+c_\Gamma(u^{n+1}-u_h^n,\theta_\phi^{n+1})
+R_u^{n+1}(\theta_u^{n+1})+R_\phi^{n+1}(\theta_\phi^{n+1})
\nonumber\\
&\quad-(\frac{\sigma_u^{n+1}-\sigma_u^{n}}{\Delta t},\theta_u^{n+1})_f
-gS_0(\frac{\sigma_\phi^{n+1}-\sigma_\phi^n}{\Delta t},\theta_\phi^{n+1})_p
-\eta\sum_{i=1}^{d-1}\int_\Gamma(\sigma_u^{n+1}\cdot\tau_i)(\theta_u^{n+1}\cdot\tau_i)ds
\nonumber\\
&\quad
+\frac{e_r^{n+1}}{exp\ttt)}(
-c_\Gamma(\theta_u^{n+1},\phi_h^n)
+c_\Gamma(u_h^n,\theta_\phi^{n+1})-a_N(u_h^n,u_h^n,\theta_u^{n+1})).  \nonumber
\end{flalign}\\
\par\textbf{Step 2. Bounds on the right-hand terms of the error equation}

Let $\varepsilon,~\varepsilon_1,~\varepsilon_2,~\varepsilon_3,~\varepsilon_r$ be several positive costants, which are independent of $h$ and $\Delta t$. Using the trace theorem, Lemma \ref{lemma 1} and the Cauchy-Schwarz inequality, we have
\begin{flalign}\label{err 2dt}
&\frac{1}{2\Delta t}\|\theta_u^{n+1}\|_{0,f}^2
-\frac{1}{2\Delta t}\|\theta_u^{n}\|_{0,f}^2
+\frac{1}{2\Delta t}\|\theta_u^{n+1}-\theta_u^{n}\|_{0,f}^2
+\nu\|\nabla\theta_u^{n+1}\|_{0,f}^2
+\frac{\eta}{2}\sum_{i=1}^{d-1}\int_\Gamma(\theta_u^{n+1}\cdot\tau_i)^2ds\nonumber\\
&
+\frac{gS_0}{2\Delta t}\|\theta_\phi^{n+1}\|_{0,p}^2
-\frac{gS_0}{2\Delta t}\|\theta_\phi^{n}\|_{0,p}^2
+\frac{gS_0}{2\Delta t}\|\theta_\phi^{n+1}-\theta_\phi^{n}\|_{0,p}^2
+gk\|\nabla\theta_\phi^{n+1}\|_{0,p}^2
\nonumber\\
&\leq 5\nu\varepsilon_1\|\nabla\theta_u^{n+1}\|_{0,f}
+2\nu\varepsilon_1\|\nabla\theta_u^{n}\|_{0,f}
+3gk\varepsilon_2\|\nabla\theta_\phi^{n+1}\|_{0,p}
+gk\varepsilon_2\|\nabla\theta_\phi^{n}\|_{0,p}
+C\Delta t\int_{t^n}^{t^{n+1}}\|\partial_tu\|_{0,f}^2dt\nonumber\\
&\quad
+C\|\theta_u^n\|_{0,f}^2+C\|\sigma_u^n\|_{0,f}^2+C\Delta t\int_{t^n}^{t^{n+1}}\|\partial_tu\|_{2,f}^2dt
+C\|\theta_u^{n+1}\|_{0,f}^2+C\|\nabla\sigma_u^n\|_{0,f}^2
\nonumber\\
&\quad+C\|\theta_u^n\|_{0,f}^2\|\nabla\theta_u^n\|_{0,f}^2
+C\|\nabla\sigma_u^n\|_{0,f}^4
+C\Delta t\int_{t^n}^{t^{n+1}}\|\partial_t\phi\|_{0,p}^2dt
+C\|\theta_\phi^n\|_{0,p}^2
+C\|\sigma_\phi^n\|_{0,p}^2\nonumber\\
&\quad
+C\Delta t\int_{t^n}^{t^{n+1}}\|\partial_t\phi\|_{1,p}^2dt
+C\|\nabla\sigma_\phi^n\|_{0,p}^2
+C\Delta t\int_{t^n}^{t^{n+1}}\|\partial_tu\|_{1,f}^2dt
+C\Delta t\int_{t^n}^{t^{n+1}}\|\partial_{tt}u\|_{0,f}^2dt
\nonumber\\
&\quad+C\Delta t\int_{t^n}^{t^{n+1}}\|\partial_{tt}\phi\|_{0,p}^2dt
+\frac{C}{\Delta t}\int_{t^n}^{t^{n+1}}\|\partial_{t}\sigma_u\|_{0,f}^2dt
+\frac{C}{\Delta t}\int_{t^n}^{t^{n+1}}\|\partial_{t}\sigma_\phi\|_{0,p}^2dt
+C\|\nabla\sigma_u^{n+1}\|_{0,f}^2\nonumber\\
&\quad
+\frac{e_r^{n+1}}{exp\ttt)}(
{-c_\Gamma(\theta_u^{n+1},\phi_h^n)
+c_\Gamma(u_h^n,\theta_\phi^{n+1})}
-a_N(u_h^n,u_h^n,\theta_u^{n+1})).&
\end{flalign}

%
See Appendix \uppercase\expandafter{\romannumeral1} for more details of the proof. The last term on the right hand side of the inequality above is difficult to bound. Therefore, to find a good way to estimate this term, we focus on the equation (\ref{er3}). \\

\textbf{Step 3. Error analysis of the auxiliary variable}

We first need to rewrite some of the trilinear terms and interfacial terms in (\ref{er3}) by
\begin{flalign*}
&{\bm {Q_1}}=\frac{e_r^{n+1}}{exp\ttt)}(a_N(u^{n+1},u^{n+1},u^{n+1})-a_N(u_h^n,u_h^n,u_h^{n+1}))\\
&\quad=\frac{e_r^{n+1}}{exp\ttt)}(a_N(u^{n+1}-u_h^n,u^{n+1},u^{n+1})+a_N(u_h^n,u^{n+1},u^{n+1})
-a_N(u_h^n,u_h^n,u^{n+1})\\
&\quad\quad+a_N(u_h^n,u_h^n,u^{n+1})-a_N(u_h^n,u_h^n,u_h^{n+1}))\\
&\quad=\frac{e_r^{n+1}}{exp\ttt)}(a_N(u^{n+1}-u_h^n,u^{n+1},u^{n+1})+a_N(u_h^n,u^{n+1}-u_h^n,u^{n+1})
\\
&\quad\quad+a_N(u_h^n,u_h^n,\theta_u^{n+1})+a_N(u_h^n,u_h^n,\sigma_u^{n+1})),&
\end{flalign*}

\begin{flalign*}
&{\bm{Q_2}}=\frac{e_r^{n+1}}{exp\ttt)}(c_\Gamma(u^{n+1},\phi^{n+1})
-c_\Gamma(u_h^{n+1},\phi_h^n))\\
&\quad=\frac{e_r^{n+1}}{exp\ttt)}(c_\Gamma(u^{n+1},\phi^{n+1}-\phi_h^n)+c_\Gamma(u^{n+1},\phi_h^n)
-c_\Gamma(u_h^{n+1},\phi_h^n))\\
&\quad=\frac{e_r^{n+1}}{exp\ttt)}(c_\Gamma(u^{n+1},\phi^{n+1}-\phi_h^n)
+c_\Gamma(\theta_u^{n+1},\phi_h^n)+c_\Gamma(\sigma_u^{n+1},\phi_h^n)),&
\end{flalign*}
\begin{flalign*}
&{\bm{Q_3}}=-\frac{e_r^{n+1}}{exp\ttt)}(c_\Gamma(u^{n+1},\phi^{n+1})
-c_\Gamma(u_h^n,\phi_h^{n+1}))\\
&\quad=-\frac{e_r^{n+1}}{exp\ttt)}(c_\Gamma(u^{n+1}-u_h^n,\phi^{n+1})+c_\Gamma(u_h^n,\phi^{n+1})
-c_\Gamma(u_h^n,\phi_h^{n+1}))\\
&\quad=-\frac{e_r^{n+1}}{exp\ttt)}(c_\Gamma(u^{n+1}-u_h^n,\phi^{n+1})+c_\Gamma(u_h^n,\theta_\phi^{n+1})
+c_\Gamma(u_h^n,\sigma_\phi^{n+1})).&
\end{flalign*}
Thus \eqref{er3} can be recast as
\begin{flalign}\label{err =}
&\frac{1}{2\Delta t}|e_r^{n+1}|^2
-\frac{1}{2\Delta t}|e_r^{n}|^2
+\frac{1}{2\Delta t}|e_r^{n+1}-e_r^{n}|^2
+\frac{1}{T}|e_r^{n+1}|^2\nonumber\\
&
=R_r^{n+1} e_r^{n+1}
+\frac{e_r^{n+1}}{exp\ttt)}(a_N(u^{n+1}-u_h^n,u^{n+1},u^{n+1})
+a_N(u_h^n,u^{n+1}-u_h^n,u^{n+1})\nonumber\\
&\quad
+{a_N(u_h^n,u_h^n,\sigma_u^{n+1})}
+c_\Gamma(u^{n+1},\phi^{n+1}-\phi_h^n)
+c_\Gamma(\sigma_u^{n+1},\phi_h^n)
-c_\Gamma(u^{n+1}-u_h^n,\phi^{n+1})\nonumber\\
&\quad
-c_\Gamma(u_h^n,\sigma_\phi^{n+1})
+c_\Gamma(\theta_u^{n+1},\phi_h^n)
-c_\Gamma(u_h^n,\theta_\phi^{n+1})
+a_N(u_h^n,u_h^n,\theta_u^{n+1})
)
.&
\end{flalign}
It is easy to see that the last three terms at the right end of the equation above will match exactly the last three terms at the right end of (\ref{err 2dt}).
Next, let us estimate the remaining terms at the right hand side of (\ref{err =}).
\begin{flalign}\label{W1}
&\bm{W_1}=R_r^{n+1} e_r^{n+1}
=(\frac{r^{n+1}-r^n}{\Delta t}-r_t^{n+1})e_r^{n+1}
\leq \left|\frac{r^{n+1}-r^n}{\Delta t}-r_t^{n+1}\right||e_r^{n+1}|\nonumber\\
&\qquad\leq \varepsilon_r|e_r^{n+1}|^2
+C\Delta t\int_{t^n}^{t^{n+1}}|r_{tt}|^2dt,&
\end{flalign}
\begin{flalign}\label{W2}
&\bm{W_2}
=\frac{e_r^{n+1}}{exp\ttt)}a_N(u^{n+1}-u_h^n,u^{n+1},u^{n+1})\nonumber\\
&\qquad\leq C|e_r^{n+1}|\|u^{n+1}-u_h^n\|_{0,f}\|u^{n+1}\|_{1,f}\|u^{n+1}\|_{2,f}\nonumber\\
&\qquad\leq \varepsilon_r|e_r^{n+1}|^2
+C\|\theta_u^n\|_{0,f}^2
+C\|\sigma_u^n\|_{0,f}^2
+C\Delta t\int_{t^n}^{t^{n+1}}\|\partial_t u\|_{0,f}^2dt,&
\end{flalign}
\begin{flalign}\label{W3}
&
{\bm{W_{3}}}
=\frac{e_r^{n+1}}{exp\ttt)}a_N(u_h^n,u^{n+1}-u_h^n,u^{n+1})\nonumber\\
&\qquad\leq C|e_r^{n+1}|\|\nabla u_h^n\|_{0,f}\|u^{n+1}-u_h^n\|_{0,f}
\|u^{n+1}\|_{2,f}\nonumber\\
&\qquad\leq
C\Delta t\int_{t^n}^{t^{n+1}}\|\partial_t u\|_{0,f}^2dt
+C\|\theta_u^n\|_{0,f}^2
+C\|\sigma_u^n\|_{0,f}^2
+\frac{1}{4C^*}\|\nabla u_h^n\|_{0,f}^2|e_r^{n+1}|^2, &
\end{flalign}
\begin{flalign}\label{W4}
&\bm{W_{4}}
=\frac{e_r^{n+1}}{exp\ttt)}a_N(u_h^n,u_h^n,\sigma_u^{n+1})\nonumber\\
&\qquad=\frac{e_r^{n+1}}{exp\ttt)}(a_N(u_h^n,u^n,\sigma_u^{n+1})
-a_N(u_h^n,\theta_u^n,\sigma_u^{n+1})-a_N(u_h^n,\sigma_u^n,\sigma_u^{n+1})).&
\end{flalign}
Using Lemma \ref{lemma 1}, the first term on the right side of (\ref{W4}) can be bounded as
\begin{flalign}\label{W4-1}
&\frac{e_r^{n+1}}{exp\ttt)}a_N(u_h^n,u^n,\sigma_u^{n+1})
\leq C|e_r^{n+1}|\|u_h^n\|_{0,f}\|u^n\|_{2,f}\|\nabla\sigma_u^{n+1}\|_{0,f}
\leq C\|\nabla\sigma_u^{n+1}\|_{0,f}^2+\varepsilon_r|e_r^{n+1}|^2.&
\end{flalign}
We can bounded the second term on the right side of (\ref{W4}) by
\begin{flalign}\label{W4-2}
&\frac{e_r^{n+1}}{exp\ttt)}a_N(u_h^n,\theta_u^n,\sigma_u^{n+1})
\leq C|e_r^{n+1}|\|\nabla u_h^n\|_{0,f}\|\nabla \theta_u^n\|_{0,f}\|\nabla\sigma_u^{n+1}\|_{0,f}\nonumber\\
&\qquad\qquad\qquad\qquad\qquad\qquad~
\leq \varepsilon_1\nu\|\nabla\theta_u^n\|_{0,f}^2
+
{C\|\nabla u_h^n\|_{0,f}^2\|\nabla \sigma_u^{n+1}\|_{0,f}^2}.&
\end{flalign}
The last term on the right side of (\ref{W4}) can be estimated by
\begin{flalign}\label{W4-3}
&\frac{e_r^{n+1}}{exp\ttt)}a_N(u_h^n,\sigma_u^n,\sigma_u^{n+1})
\leq C|e_r^{n+1}|\|\nabla u_h^n\|_{0,f}
\|\nabla\sigma_u^n\|_{0,f}
\|\nabla\sigma_u^{n+1}\|_{0,f}\nonumber\\
&\qquad\qquad\qquad\qquad\qquad\qquad~
\leq C\|\nabla\sigma_u^{n+1}\|_{0,f}^2
+
{C\|\nabla u_h^n\|_{0,f}^2\|\nabla\sigma_u^n\|_{0,f}^2}.&
\end{flalign}
In addition,
\begin{flalign}\label{W5}
&\bm{W_{5}}=\frac{e_r^{n+1}}{exp\ttt)}
(c_\Gamma(u^{n+1},\phi^{n+1}-\phi_h^n)+c_\Gamma(\sigma_u^{n+1},\phi_h^n)
-c_\Gamma(u^{n+1}-u_h^n,\phi^{n+1})-c_\Gamma(u_h^n,\sigma_\phi^{n+1}))\nonumber\\
&\qquad
=\frac{e_r^{n+1}}{exp\ttt)}
(c_\Gamma(u^{n+1},\phi^{n+1}-\phi_h^n)
+c_\Gamma(\sigma_u^{n+1},\phi^n-\theta_\phi^n-\sigma_\phi^n)
-c_\Gamma(u^{n+1}-u_h^n,\phi^{n+1})\nonumber\\
&\qquad\quad-c_\Gamma(u^n-\theta_u^n-\sigma_u^n,\sigma_\phi^{n+1}))\nonumber\\
&\qquad\leq
\varepsilon_r|e_r^{n+1}|^2
+C\|\phi^{n+1}-\phi_h^n\|_{0,p}\|\nabla(\phi^{n+1}-\phi_h^n)\|_{0,p}
+C\|\nabla\sigma_u^{n+1}\|_{0,f}^2
+\varepsilon_2gk\|\nabla\theta_\phi^n\|_{0,p}^2
+C\|\nabla\sigma_\phi^n\|_{0,p}^2\nonumber\\
&\qquad\quad+C\|u^{n+1}-u_h^n\|_{0,f}\|\nabla(u^{n+1}-u_h^n)\|_{0,f}
+C\|\nabla\sigma_\phi^{n+1}\|_{0,p}^2
+\varepsilon_1\nu\|\nabla\theta_u^n\|_{0,f}^2
+C\|\nabla\sigma_u^n\|_{0,f}^2
\nonumber\\
&\qquad\leq \varepsilon_r|e_r^{n+1}|^2
+C\Delta t\int_{t^n}^{t^{n+1}}\|\partial_t\phi\|_{0,p}^2dt
+C\|\theta_\phi^n\|_{0,p}^2+C\|\sigma_\phi^n\|_{0,p}^2
+C\Delta t\int_{t^n}^{t^{n+1}}\|\partial_t\phi\|_{1,p}^2dt\nonumber\\
&\qquad\quad
+\varepsilon_2gk\|\nabla\theta_\phi^n\|_{0,p}^2
+C\|\nabla\sigma_\phi^n\|_{0,p}^2+C\|\nabla\sigma_u^{n+1}\|_{0,f}^2
+C\|\nabla\sigma_\phi^n\|_{0,p}^2
+C\Delta t\int_{t^n}^{t^{n+1}}\|\partial_tu\|_{0,f}^2dt
\nonumber\\
&\qquad\quad+C\|\theta_u^n\|_{0,f}^2
+C\|\sigma_u^n\|_{0,f}^2
+C\Delta t\int_{t^n}^{t^{n+1}}\|\partial_tu\|_{1,f}^2dt
+\varepsilon_1\nu\|\nabla\theta_u^n\|_{0,f}^2
+C\|\nabla\sigma_\phi^{n+1}\|_{0,p}^2
+C\|\nabla\sigma_u^n\|_{0,f}^2.&
\end{flalign}
Choosing $\varepsilon_r=\frac{1}{8T},\varepsilon_1=\frac{1}{18},\varepsilon_2=\frac{1}{10}$ and combining \eqref{err 2dt} with (\ref{err =})-(\ref{W5}), we have
\begin{flalign}\label{toget}
&\frac{1}{2\Delta t}\|\theta_u^{n+1}\|_{0,f}^2
-\frac{1}{2\Delta t}\|\theta_u^{n}\|_{0,f}^2
+\frac{1}{2\Delta t}\|\theta_u^{n+1}-\theta_u^{n}\|_{0,f}^2
+\frac{\nu}{2}\|\nabla\theta_u^{n+1}\|_{0,f}^2
+\frac{2\nu}{9}\left(\|\nabla\theta_u^{n+1}\|_{0,f}^2
-\|\nabla\theta_u^{n}\|_{0,f}^2\right)\nonumber\\
&+\frac{\eta}{2}\sum_{i=1}^{d-1}\int_\Gamma(\theta_u^{n+1}\cdot\tau_i)^2ds
+\frac{gS_0}{2\Delta t}\|\theta_\phi^{n+1}\|_{0,p}^2
-\frac{gS_0}{2\Delta t}\|\theta_\phi^{n}\|_{0,p}^2
+\frac{gS_0}{2\Delta t}\|\theta_\phi^{n+1}-\theta_\phi^{n}\|_{0,p}^2
+\frac{gk}{2}\|\nabla\theta_\phi^{n+1}\|_{0,p}^2\nonumber\\
&+\frac{gk}{5}(\|\nabla\theta_\phi^{n+1}\|_{0,p}^2
-\|\nabla\theta_\phi^{n}\|_{0,p}^2)
+\frac{1}{2\Delta t}|e_r^{n+1}|^2
-\frac{1}{2\Delta t}|e_r^{n}|^2
+\frac{1}{2\Delta t}|e_r^{n+1}-e_r^{n}|^2
+\frac{1}{2T}|e_r^{n+1}|^2\nonumber\\
&\leq
C\Delta t\int_{t^n}^{t^{n+1}}\|\partial_tu\|_{0,f}^2dt
+C\|\theta_u^n\|_{0,f}^2+C\|\sigma_u^n\|_{0,f}^2
+C\Delta t\int_{t^n}^{t^{n+1}}\|\partial_tu\|_{2,f}^2dt
+C\|\theta_u^{n+1}\|_{0,f}^2\nonumber\\
&\quad
+C\|\nabla\sigma_u^n\|_{0,f}^2
+C\|\theta_u^n\|_{0,f}^2\|\nabla\theta_u^n\|_{0,f}^2
+C\|\nabla\sigma_u^n\|_{0,f}^4
+C\Delta t\int_{t^n}^{t^{n+1}}\|\partial_t\phi\|_{0,p}^2dt
+C\|\theta_\phi^n\|_{0,p}^2\nonumber\\
&\quad+C\|\sigma_\phi^n\|_{0,p}^2
+C\Delta t\int_{t^n}^{t^{n+1}}\|\partial_t\phi\|_{1,p}^2dt
+C\|\nabla\sigma_\phi^n\|_{0,p}^2
+C\Delta t\int_{t^n}^{t^{n+1}}\|\partial_tu\|_{1,f}^2dt
\nonumber\\
&\quad+C\Delta t\int_{t^n}^{t^{n+1}}\|\partial_{tt}u\|_{0,f}^2dt
+C\Delta t\int_{t^n}^{t^{n+1}}\|\partial_{tt}\phi\|_{0,p}^2dt
+\frac{C}{\Delta t}\int_{t^n}^{t^{n+1}}\|\partial_{t}\sigma_u\|_{0,f}^2dt
\nonumber\\
&\quad
+\frac{C}{\Delta t}\int_{t^n}^{t^{n+1}}\|\partial_{t}\sigma_\phi\|_{0,p}^2dt
+C\|\nabla\sigma_u^{n+1}\|_{0,f}^2
+C\Delta t\int_{t^n}^{t^{n+1}}|r_{tt}|^2dt
+C\|\nabla u_h^n\|_{0,f}^2\|\nabla \sigma_u^{n+1}\|_{0,f}^2\nonumber\\
&\quad
+C\|\nabla u_h^n\|_{0,f}^2\|\nabla \sigma_u^{n}\|_{0,f}^2
+C\|\nabla \sigma _\phi^{n+1}\|_{0,p}^2
+\frac{1}{4C^*}\|\nabla u_h^n\|_{0,f}^2|e_r^{n+1}|^2.&
\end{flalign}
Let $|e_r^M|^2$ denote the maximum value of $|e_r^{n+1}|^2~(n=0,\cdots,N-1)$. Then summing up (\ref{toget}) from $n=0$ to $n=M-1$ yields
\begin{flalign}\label{err * 2dt}
&\|\theta_u^M\|_{0,f}^2
+\sum_{n=0}^{M-1}\|\theta_u^{n+1}-\theta_u^n\|_{0,f}^2
+\nu\Delta t\sum_{n=0}^{M-1}\|\nabla\theta_u^{n+1}\|_{0,f}^2
+\frac{4\nu\Delta t}{9}\|\nabla\theta_u^M\|_{0,f}^2
+\frac{\eta}{2}\Delta t\sum_{i=1}^{d-1}\int_\Gamma(\theta_u^M\cdot\tau_i)^2ds\nonumber\\
&+gS_0\|\theta_\phi^M\|_{0,p}^2
+gS_0\sum_{n=0}^{M-1}\|\theta_\phi^{n+1}-\theta_\phi^n\|_{0,p}^2
+gk\Delta t\sum_{n=0}^{M-1}\|\nabla\theta_\phi^{n+1}\|_{0,p}^2
+\frac{2gk}{5}\Delta t\|\nabla\theta_\phi^M\|_{0,p}^2
+|e_r^M|^2\nonumber\\
&+\sum_{n=0}^{M-1}|e_r^{n+1}-e_r^n|^2
+\frac{\Delta t}{T}\sum_{n=0}^{M-1}|e_r^{n+1}|^2\nonumber\\
&\leq
\sum_{n=0}^{M-1}\bigg\{
C\Delta t^2\int_{t^n}^{t^{n+1}}\|\partial_tu\|_{0,f}^2dt
+C\Delta t\|\theta_u^n\|_{0,f}^2
+C\Delta t\|\sigma_u^n\|_{0,f}^2
+C\Delta t^2\int_{t^n}^{t^{n+1}}\|\partial_tu\|_{2,f}^2dt
\nonumber\\
&\quad
+C\Delta t\|\theta_u^{n+1}\|_{0,f}^2
+C\Delta t\|\nabla\sigma_u^{n}\|_{0,f}^2
+C\Delta t\|\theta_u^{n}\|_{0,f}^2\|\nabla\theta_u^{n}\|_{0,f}^2
+C\Delta t\|\nabla\sigma_u^n\|_{0,f}^4
\nonumber\\
&\quad
+C\Delta t^2\int_{t^n}^{t^{n+1}}\|\partial_t\phi\|_{0,p}^2dt
+C\Delta t\|\theta_\phi^n\|_{0,p}^2
+C\Delta t\|\sigma_\phi^n\|_{0,p}^2
+C\Delta t^2\int_{t^n}^{t^{n+1}}\|\partial_t\phi\|_{1,p}^2dt
\nonumber\\
&\quad+C\Delta t\|\nabla\sigma_\phi^n\|_{0,p}^2
+C\Delta t^2\int_{t^n}^{t^{n+1}}\|\partial_tu\|_{1,f}^2dt
+C\Delta t^2\int_{t^n}^{t^{n+1}}\|\partial_{tt}u\|_{0,f}^2dt
+C\Delta t^2\int_{t^n}^{t^{n+1}}\|\partial_{tt}\phi\|_{0,p}^2dt\nonumber\\
&\quad
+C\int_{t^n}^{t^{n+1}}\|\partial_{t}\sigma_u\|_{0,f}^2dt
+C\int_{t^n}^{t^{n+1}}\|\partial_{t}\sigma_\phi\|_{0,p}^2dt
+C\Delta t\|\nabla\sigma_u^{n+1}\|_{0,f}^2
+C\Delta t^2\int_{t^n}^{t^{n+1}}|r_{tt}|^2dt
\nonumber\\
&\quad
+C\Delta t\|\nabla u_h^n\|_{0,f}^2\|\nabla\sigma_u^{n+1}\|_{0,f}^2
+C\Delta t\|\nabla u_h^n\|_{0,f}^2\|\nabla\sigma_u^{n}\|_{0,f}^2
+C\Delta t\|\nabla \sigma_\phi^{n+1}\|_{0,p}^2
\bigg\}\nonumber\\
&\quad
+\Delta t\frac{1}{2C^*}\sum_{n=0}^{M-1}\|\nabla u_h^n\|_{0,f}^2|e_r^M|^2.
&
\end{flalign}
Recalling the estimates of Ritz projection operators (\ref{Ritz u}) and (\ref{Ritz phi}), we derive that
\begin{flalign}\label{before gronwall}
&\|\theta_u^M\|_{0,f}^2
+\sum_{n=0}^{M-1}\|\theta_u^{n+1}-\theta_u^n\|_{0,f}^2
+\nu\Delta t\sum_{n=0}^{M-1}\|\nabla\theta_u^{n+1}\|_{0,f}^2
+\frac{4\nu\Delta t}{9}\|\nabla\theta_u^M\|_{0,f}^2
+\frac{\eta}{2}\Delta t\sum_{i=1}^{d-1}\int_\Gamma(\theta_u^M\cdot\tau_i)^2ds\nonumber\\
&+gS_0\|\theta_\phi^M\|_{0,p}^2
+gS_0\sum_{n=0}^{M-1}\|\theta_\phi^{n+1}-\theta_\phi^n\|_{0,p}^2
+gk\Delta t\sum_{n=0}^{M-1}\|\nabla\theta_\phi^{n+1}\|_{0,p}^2
+\frac{2gk}{5}\Delta t\|\nabla\theta_\phi^M\|_{0,p}^2
+\frac{1}{2}|e_r^M|^2\nonumber\\
&+\sum_{n=0}^{M-1}|e_r^{n+1}-e_r^n|^2
+\frac{\Delta t}{T}\sum_{n=0}^{M-1}|e_r^{n+1}|^2\nonumber\\
&
\leq \Delta t\sum_{n=0}^{N}C(\|\theta_u^n\|_{0,f}^2+\|\theta_\phi^n\|_{0,p}^2)
+\Delta t\sum_{n=0}^{N}C\|\nabla\theta_u^{n}\|_{0,f}^2\|\theta_u^{n}\|_{0,f}^2\nonumber\\
&\quad
+C\Delta t^2
(1+\|\partial_t u\|^2_{2,0,f}+\|\partial_t u\|^2_{2,1,f}+\|\partial_t u\|^2_{2,2,f}
+\|\partial_{tt} u\|^2_{2,0,f}
+\|\partial_t \phi\|^2_{2,0,p}
+\|\partial_t \phi\|^2_{2,1,p}\nonumber\\
&\quad+\|\partial_{tt} \phi\|^2_{2,0,p})
+Ch^{2k}\normmm{u}^2_{2,k+1,f}
+Ch^{2k}\normmm{u}^2_{\infty,k+1,f}
+Ch^{2k+2}\|\partial_t u\|^2_{2,k+1,f}
+Ch^{2m}\normmm{\phi}^2_{2,m+1,p}\nonumber\\
&\quad
+Ch^{2m+2}\|\partial_t \phi\|^2_{2,m+1,p}.&
\end{flalign}

\par\textbf{Step 4. Adopt the Gronwall lemma}

We can obtain the following inequality by using the Gronwall lemma:
\begin{flalign}\label{after gronwall}
&\|\theta_u^M\|_{0,f}^2
+\sum_{n=0}^{M-1}\|\theta_u^{n+1}-\theta_u^n\|_{0,f}^2
+\nu\Delta t\sum_{n=0}^{M-1}\|\nabla \theta_u^{n+1}\|_{0,f}^2
+\frac{4\nu}{9}\Delta t\|\nabla \theta_u^M\|_{0,f}^2
+\frac{\eta}{2}\Delta t\sum_{i=1}^{d-1}\int_\Gamma(\theta_u^M\cdot\tau_i)^2ds\nonumber\\
&+gS_0\|\theta_\phi^M\|_{0,p}^2
+gS_0\sum_{n=0}^{M-1}\|\theta_\phi^{n+1}-\theta_\phi^n\|_{0,p}^2
+gk\Delta t\sum_{n=0}^{M-1}\|\nabla \theta_\phi^{n+1}\|_{0,p}^2
+\frac{2gk}{5}\Delta t\|\nabla \theta_\phi^M\|_{0,p}^2
+\frac{1}{2}|e_r^M|^2\nonumber\\
&+\sum_{n=0}^{M-1}|e_r^{n+1}-e_r^n|^2
+\frac{\Delta t}{T}\sum_{n=0}^{M-1}|e_r^{n+1}|^2\nonumber\\
&\leq
C\Delta t^2
(1+\|\partial_t u\|^2_{2,0,f}+\|\partial_t u\|^2_{2,1,f}+\|\partial_t u\|^2_{2,2,f}
+\|\partial_{tt} u\|^2_{2,0,f}+\|\partial_t \phi\|^2_{2,0,p}
+\|\partial_t \phi\|^2_{2,1,p}+\|\partial_{tt} \phi\|^2_{2,0,p})\nonumber\\
&\quad
+Ch^{2k}\normmm{u}^2_{2,k+1,f}
+Ch^{2k}\normmm{u}^2_{\infty,k+1,f}
+Ch^{2k+2}\|\partial_t u\|^2_{2,k+1,f}
+Ch^{2m}\normmm{\phi}^2_{2,m+1,p}\nonumber\\
&\quad
+Ch^{2m+2}\|\partial_t \phi\|^2_{2,m+1,p}.&
\end{flalign}
Since $|e_r^M|^2=\max_{1\leq s\leq N}|e_r^s|^2$, there holds
\begin{flalign}\label{err max}
&\|\theta_u^s\|_{0,f}^2
+\sum_{n=0}^{s-1}\|\theta_u^{n+1}-\theta_u^n\|_{0,f}^2
+\nu\Delta t\sum_{n=0}^{s-1}\|\nabla \theta_u^{n+1}\|_{0,f}^2
+\frac{4\nu}{9}\Delta t\|\nabla \theta_u^s\|_{0,f}^2
+\frac{\eta}{2}\Delta t\sum_{i=1}^{d-1}\int_\Gamma(\theta_u^s\cdot\tau_i)^2ds\nonumber\\
&+gS_0\|\theta_\phi^s\|_{0,p}^2
+gS_0\sum_{n=0}^{s-1}\|\theta_\phi^{n+1}-\theta_\phi^n\|_{0,p}^2
+gk\Delta t\sum_{n=0}^{s-1}\|\nabla\theta_\phi^{n+1}\|_{0,p}^2
+\frac{2gk}{5}\Delta t\|\nabla\theta_\phi^s\|_{0,p}^2
+\frac{1}{2}|e_r^s|^2\nonumber\\
&+\sum_{n=0}^{s-1}|e_r^{n+1}-e_r^n|^2
+\frac{\Delta t}{T}\sum_{n=0}^{s-1}|e_r^{n+1}|^2\nonumber\\
&\leq
C\Delta t^2
(1+\|\partial_t u\|^2_{2,0,f}+\|\partial_t u\|^2_{2,1,f}+\|\partial_t u\|^2_{2,2,f}
+\|\partial_{tt} u\|^2_{2,0,f}+\|\partial_t \phi\|^2_{2,0,p}
+\|\partial_t \phi\|^2_{2,1,p}+\|\partial_{tt} \phi\|^2_{2,0,p})\nonumber\\
&\quad
+Ch^{2k}\normmm{u}^2_{2,k+1,f}
+Ch^{2k}\normmm{u}^2_{\infty,k+1,f}
+Ch^{2k+2}\|\partial_t u\|^2_{2,k+1,f}
+Ch^{2m}\normmm{\phi}^2_{2,m+1,p}\nonumber\\
&\quad
+Ch^{2m+2}\|\partial_t \phi\|^2_{2,m+1,p}.&
\end{flalign}
Also, it is easy to verify that
\begin{flalign}\label{sigma}
&\|\sigma_u^s\|_{0,f}^2
+\nu\Delta t\sum_{n=0}^{s-1}\|\nabla\sigma_u^{n+1}\|_{0,f}^2
+\frac{4\nu}{9}\Delta t\|\nabla\sigma_u^s\|_{0,f}^2
+\frac{\eta}{2}\Delta t\sum_{i=1}^{d-1}\int_\Gamma(\sigma_u^s\cdot\tau_i)^2ds
+gS_0\|\sigma_\phi^s\|_{0,p}^2\nonumber\\
&
+\frac{2gk}{5}\Delta t\|\nabla\sigma_\phi^s\|_{0,p}^2\nonumber\\
&\leq Ch^{2k+2}\normmm{u}_{\infty,k+1,f}^2
+Ch^{2k}\normmm{u}_{2,k+1,f}^2
+Ch^{2m+2}\normmm{\phi}_{\infty,m+1,p}^2
+Ch^{2m}\normmm{\phi}_{2,m+1,p}^2,&
\end{flalign}
where the desired result \eqref{err 1order} can be obtained using the triangle inequality.
\end{proof}

\section{Numerical Tests}
In this section, we will present some numerical examples of the constructed first- and second-order schemes to verify our theoretical results and to show the efficiency for the unsteady Navier-Stokes-Darcy model.

To satisfy the discrete inf-sup condition, we select Taylor-Hood finite element (P2-P1) for the velocity and pressure for the free fluid, and choose P2 element for the hydraulic head of porous media flow, then by using Theorem \ref{final thm}, we can derive that
\begin{align}
&\|e_u^s\|_{0,f}^2
+\sum_{n=0}^{s-1}\|e_u^{n+1}-e_u^n\|_{0,f}^2
+\nu\Delta t\sum_{n=0}^{s-1}\|\nabla e_u^{n+1}\|_{0,f}^2
+\frac{4\nu}{9}\Delta t\|\nabla e_u^s\|_{0,f}^2
+\frac{\eta}{2}\Delta t\sum_{i=1}^{d-1}\int_\Gamma(e_u^s\cdot\tau_i)^2ds\nonumber\\
&+gS_0\|e_\phi^s\|_{0,p}^2
+gS_0\sum_{n=0}^{s-1}\|e_\phi^{n+1}-e_\phi^n\|_{0,p}^2
+gk\Delta t\sum_{n=0}^{s-1}\|\nabla e_\phi^{n+1}\|_{0,p}^2
+\frac{2gk}{5}\Delta t\|\nabla e_\phi^s\|_{0,p}^2
+\frac{1}{2}|e_r^s|^2\nonumber\\
&+\sum_{n=0}^{s-1}|e_r^{n+1}-e_r^n|^2
+\frac{\Delta t}{T}\sum_{n=0}^{s-1}|e_r^{n+1}|^2\nonumber\\
&\leq
Ch^4+C\Delta t^2,
\end{align}
where the positive constant $C$ is independent of $h$ and $\Delta t$.

\subsection{Convergence test}
Let $\Omega_f=[0,1]\times[0,1]$, $\Omega_p=[0,1]\times[-1,0]$, and $\Gamma=[0,1]\times\{0\}$.
We select the boundary conditions and initial solutions satisfying the following analytic solution
\begin{align*}
&u=[c sin^2(\pi x)sin(2\pi y)t^4,- c sin(2\pi x)sin^2(\pi y)t^4],\\
&p=cy cos(\pi x)t^4,\\
&\phi=csin(\pi x)sin^2(\pi y)t^4.
\end{align*}
We set $T=1$, $\alpha=1$, $S_0=1$, $g=1$, $\nu=10^{-3}$, $\mathcal{K}=k\mathbb{I}$, $k=0.1$, $c=0.01$.
In order to test the error convergence rates of the numerical schemes constructed in this paper, we take $h^2=\Delta t$ and $h=\Delta t$ for the first- and second-order schemes respectively. Errors and convergence rates  for the velocity, pressure and hydraulic head of the first- and second-order schemes are presented in Tables \ref{test1_1} and \ref{test1_2}, which are consistent with our theoretical results in Theorem \ref{final thm}.

\begin{table}[htbp]
  \centering
    \caption{Numerical solution errors and convergence rates for Algorithm 1.}\label{test1_1}
  \begin{tabular}{ccccccccc}
    \hline
    $\Delta t$ & $\|u-u_h\|_{\textit{l}^2(H^1)}$ & ${\rm Rate}$ & $\|p-p_h\|_{\textit{l}^\infty(L^2)}$ & ${\rm Rate}$ & $\|\phi-\phi_h\|_{\textit{l}^2(H^1)}$ & ${\rm Rate}$\\
    \hline
     $1/4$ & $2.452\times 10^{-2}$ & $-$ & $6.862\times 10^{-4}$ & $-$ & $5.456\times 10^{-3}$ & $-$\\
    $1/16$ & $4.693\times 10^{-3}$	&$1.193$	&$1.822\times 10^{-4}$	&$0.957$	&$1.103\times 10^{-3}$	&$1.153$\\
    $1/64$ & $1.280\times 10^{-3}$	&$0.937$	&$4.743\times 10^{-5}$	&$0.971$	&$2.596\times 10^{-4}$	&$1.044$\\
    $1/256$ & $3.054\times 10^{-4}$	&$1.034$	&$1.202\times 10^{-5}$	&$0.990$	&$6.392\times 10^{-5}$	&$1.011$\\
    $1/1024$ & $6.668\times 10^{-5}$	&$1.098$	&$3.016\times 10^{-6}$	&$0.997$	&$1.592\times 10^{-5}$	&$1.003$\\

    \hline
  \end{tabular}
\label{table1}
\end{table}

\begin{table}[h]
  \centering
    \caption{Numerical solution errors and convergence rates by Algorithm 2.}\label{test1_2}
  \begin{tabular}{ccccccccc}
    \hline
    $\Delta t$ & $\|u-u_h\|_{\textit{l}^2(H^1)}$ & ${\rm Rate}$ & $\|p-p_h\|_{\textit{l}^\infty(L^2)}$ & ${\rm Rate}$ & $\|\phi-\phi_h\|_{\textit{l}^2(H^1)}$ & ${\rm Rate}$\\
    \hline
    $1/4$ & $7.219\times 10^{-3}$	&$-$	&$1.435\times 10^{-4}$	&$-$	&$1.635\times 10^{-3}$	&$-$\\
    $1/8$ & $1.872\times 10^{-3}$	&$1.947$	&$3.359\times 10^{-5}$	&$2.095$	&$3.962\times 10^{-4}$	&$2.045$\\
    $1/16$ & $4.304\times 10^{-4}$	&$2.121$	&$8.379\times 10^{-6}$	&$2.003$	&$9.721\times 10^{-5}$	&$2.027$\\
    $1/32$ & $9.638\times 10^{-5}$	&$2.159$	&$2.103\times 10^{-6}$	&$1.995$	&$2.408\times 10^{-5}$ &$2.013$\\
    $1/64$ & $3.173\times 10^{-5}$	&$1.603$	&$5.458\times 10^{-7}$	&$1.946$	&$6.000\times 10^{-6}$	&$2.005$\\

    \hline
  \end{tabular}
\label{table2}
\end{table}

\subsection{Simulation of fluid flow in a Y-shape domain}
\begin{figure}[htbp]
  \centering
  \begin{minipage}[t]{0.49\linewidth}
  \includegraphics[width=7.1cm]{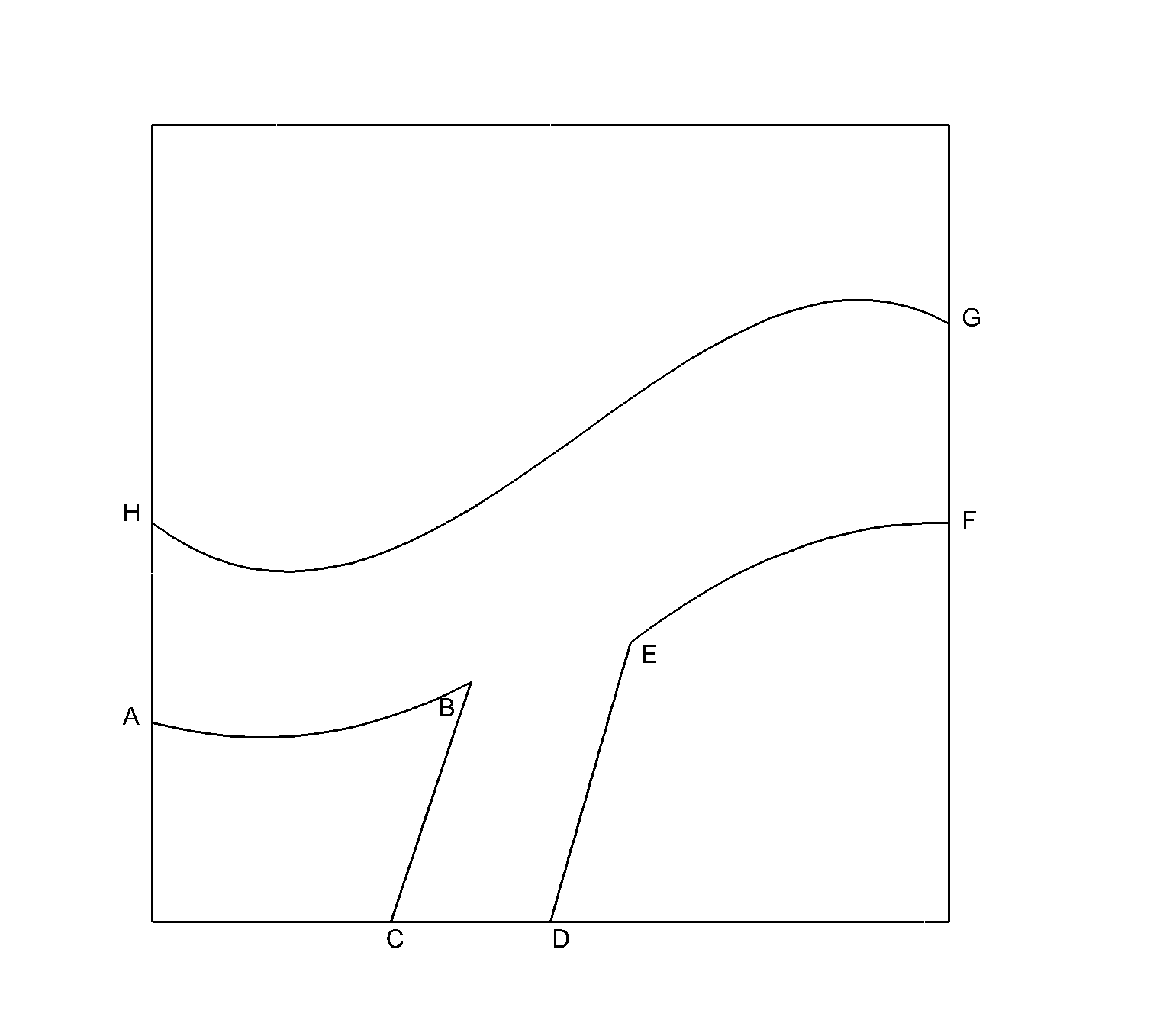}
  \end{minipage}
  \begin{minipage}[t]{0.49\linewidth}
  \includegraphics[width=7.1cm]{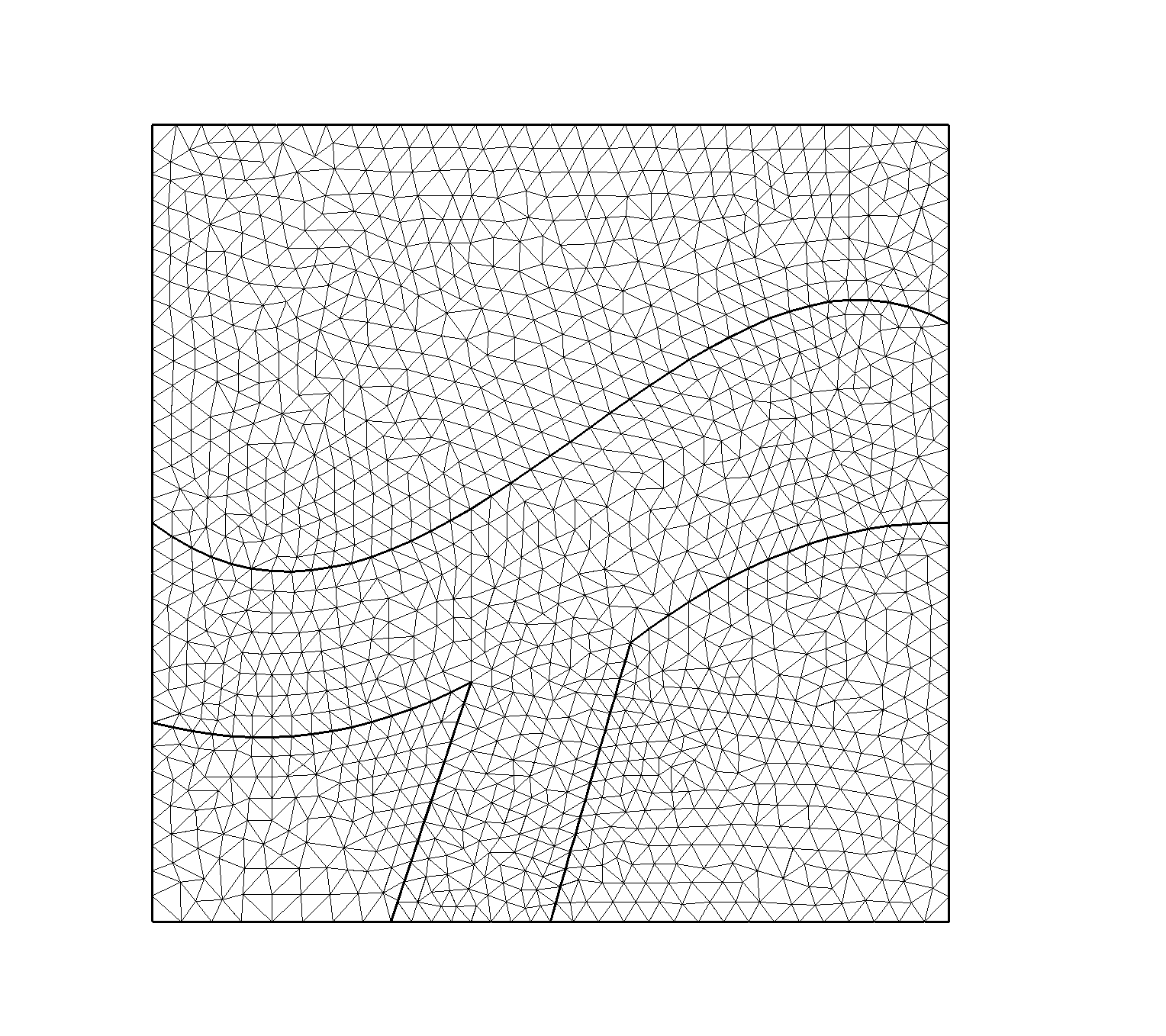}
  \end{minipage}
    \caption{Schematic diagram of the Y-shape domain and mesh partition.} \label{mesh}
\end{figure}

This section presents a simulation of subsurface flow in a karst aquifer \cite{jiang2021SAV} using Algorithm 1. The computational domain $\Omega=[0,1]\times[0,1]$ is divided into two subdomains $\Omega_f$ and $\Omega_p$. Let $\Omega_f$ be a bounded Y-shape domain $\overline{ABCDEFGH}$ with $A=(0,0.25),~B=(0.4,0.3),C=(0.3,0),D=(0.5,0),E=(0.6,0.35),
F=(1,0.5),G=(1,0.75),H=(0,0.5)$, and $\Omega_p$ be the rest parts of $\Omega$, as shown in Figure \ref{mesh}. Take $T=0.5$, $\Delta t=0.01$, $\alpha,~\nu,~g,~S_0$ all be 1. We divided the Y-shaped free flow domain into 954 triangular elements and the remaining porous media domain into 1812 triangular elements, respectively. The homogeneous Dirichlet condition is imposed on the hydraulic head $\phi=0$ on $\partial\Omega_p\backslash \Gamma$. Except for the velocity boundary conditions defined below, all values are zero.

\begin{equation*}
u=\left\{
\begin{aligned}
&(\omega_1,0),\quad {\rm on}~\overline{HA},\\
&(0,\omega_1),\quad {\rm on}~\overline{CD},\\
&(\omega_2,0),\quad {\rm on}~\overline{FG}.
\end{aligned}\right.
\end{equation*}

\begin{figure}[htbp]
  \centering
  \begin{minipage}[t]{0.49\linewidth}
  \centering
  \subfloat[][Velocity]{\label{equalVk1}
  \includegraphics[width=7.1cm]{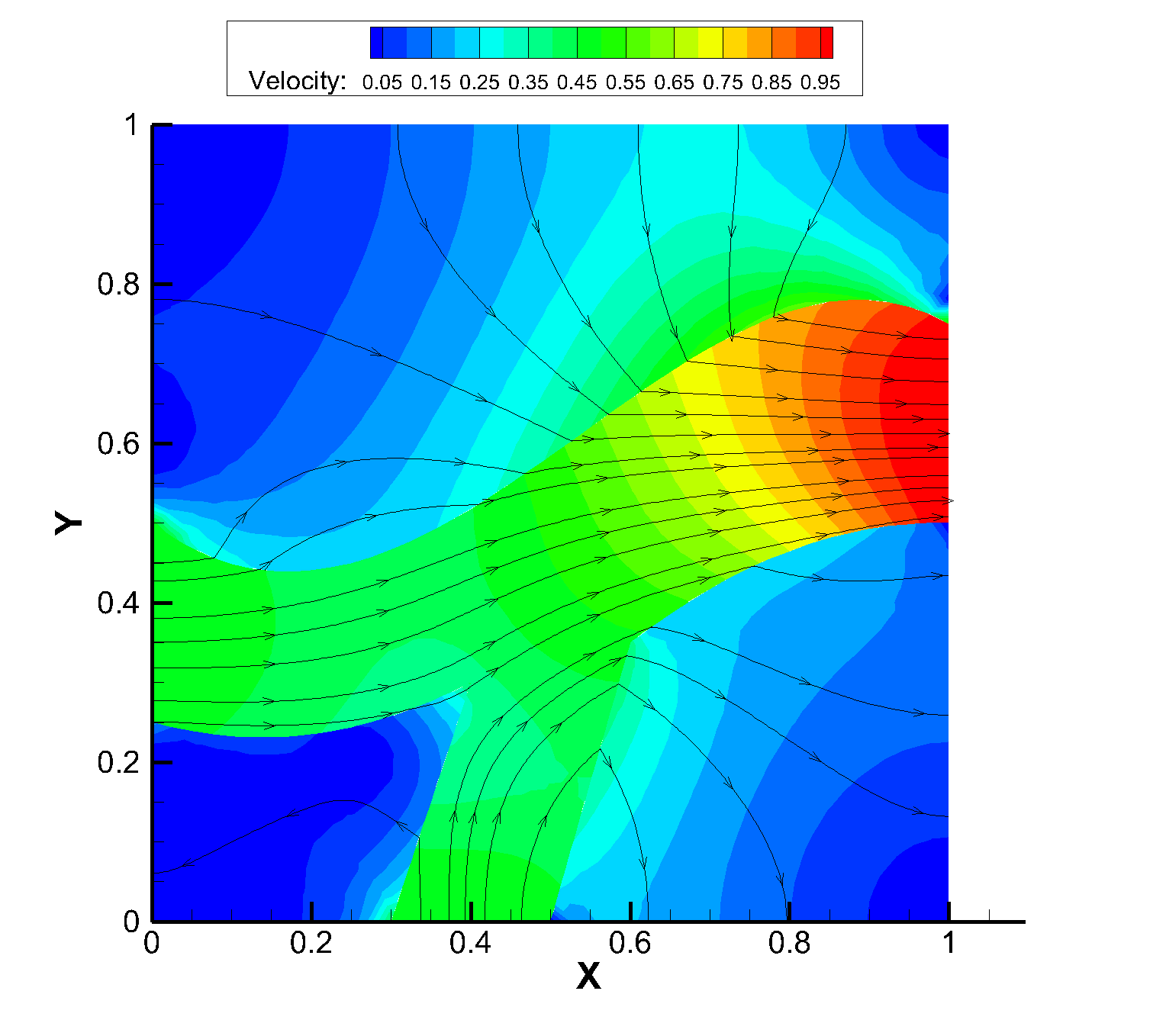}}
  \end{minipage}
  \begin{minipage}[t]{0.49\linewidth}
  \centering
  \subfloat[][Pressure]{\label{equalPk1}
  \includegraphics[width=7.1cm]{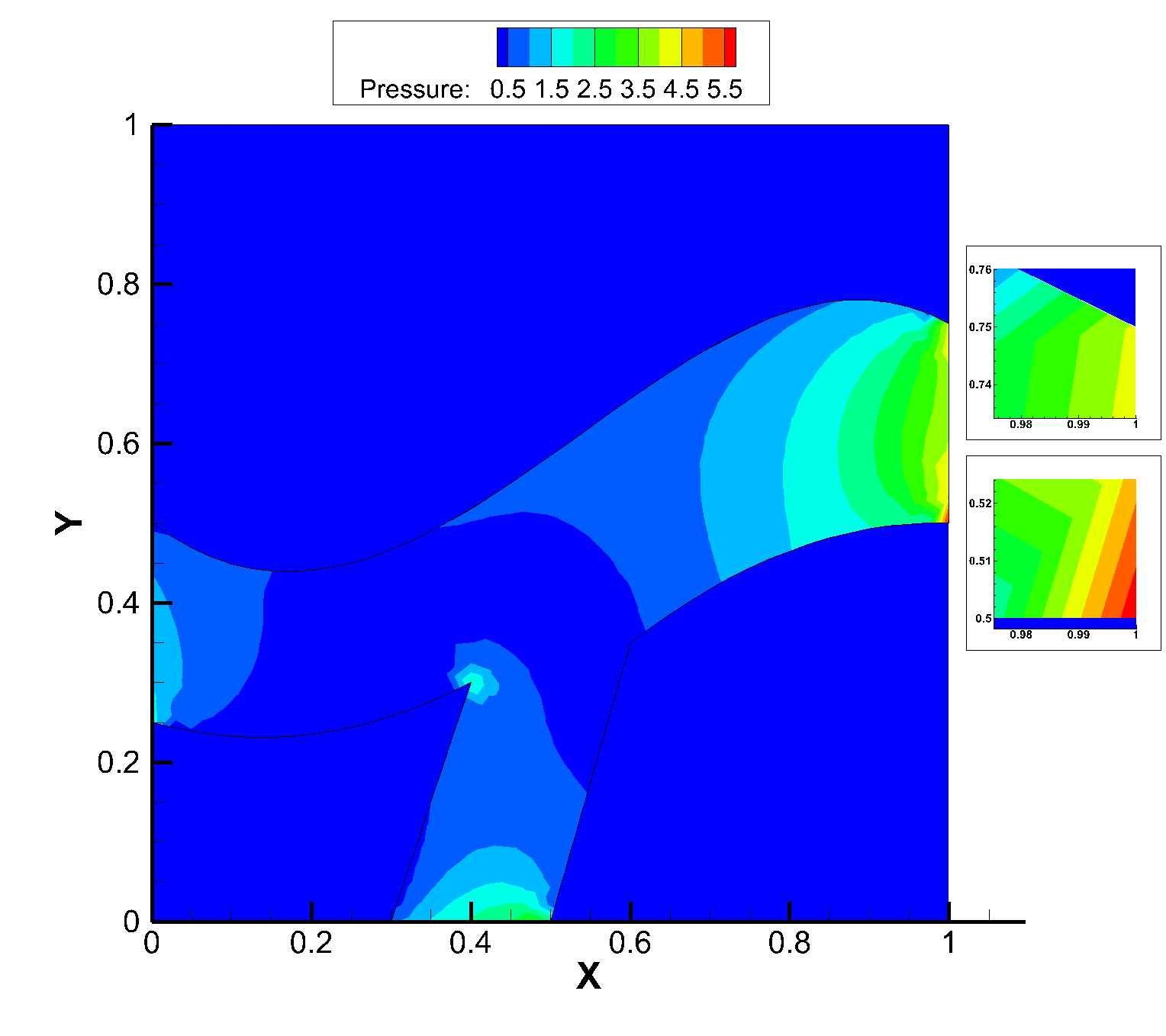}}
  \end{minipage}
  \caption{ Evolutions of velocity and pressure for Algorithm 1 ($\omega_1=0.5,~\omega_2=1,~k=1$).}\label{k=1}
\end{figure}

\begin{figure}[htbp]
  \centering
  \begin{minipage}[t]{0.49\linewidth}
  \centering
  \subfloat[][Velocity]{\label{equalVk-2}
  \includegraphics[width=7.1cm]{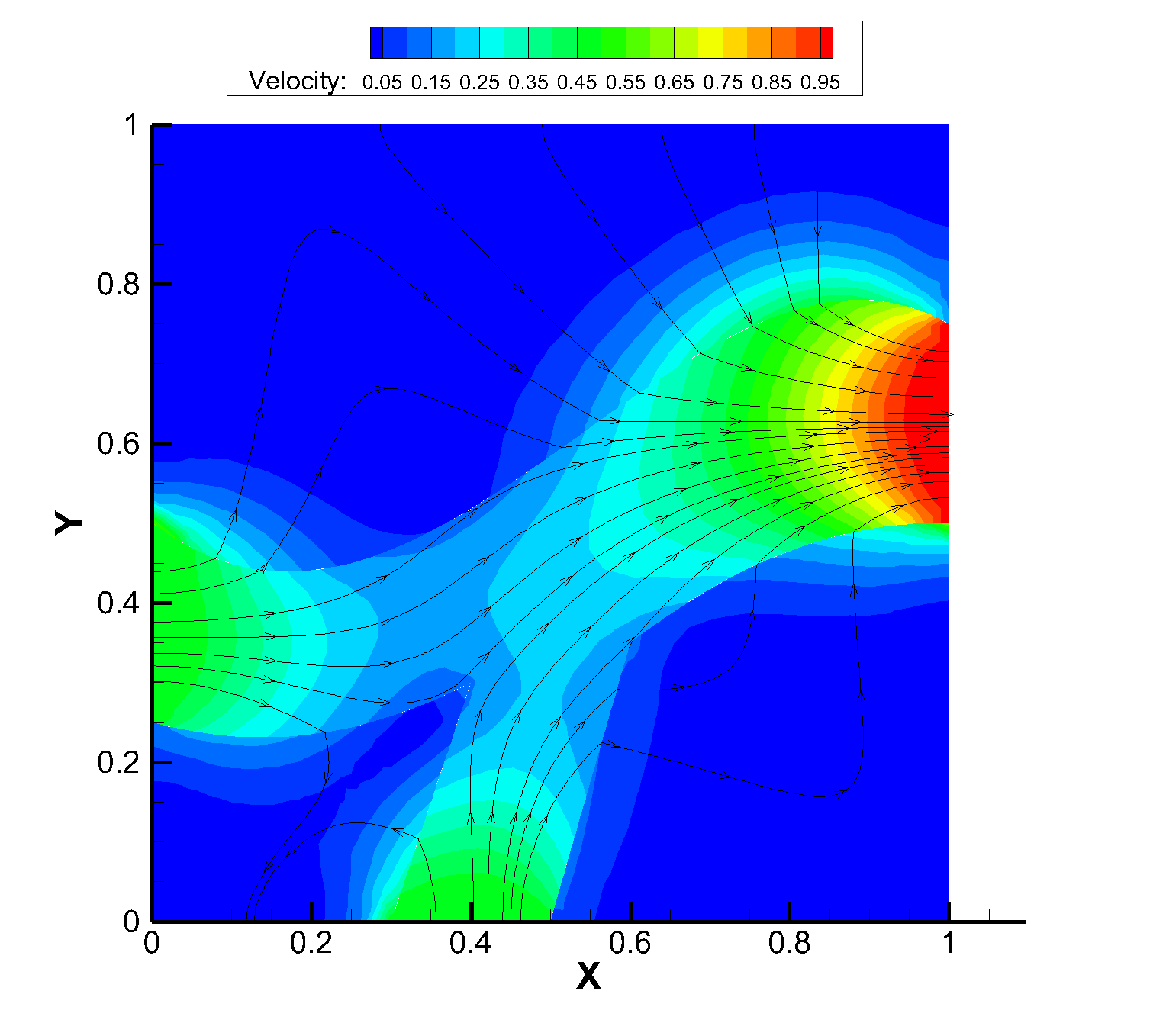}}
  \end{minipage}
  \begin{minipage}[t]{0.49\linewidth}
  \centering
  \subfloat[][Pressure]{\label{equalPk-2}
  \includegraphics[width=7.1cm]{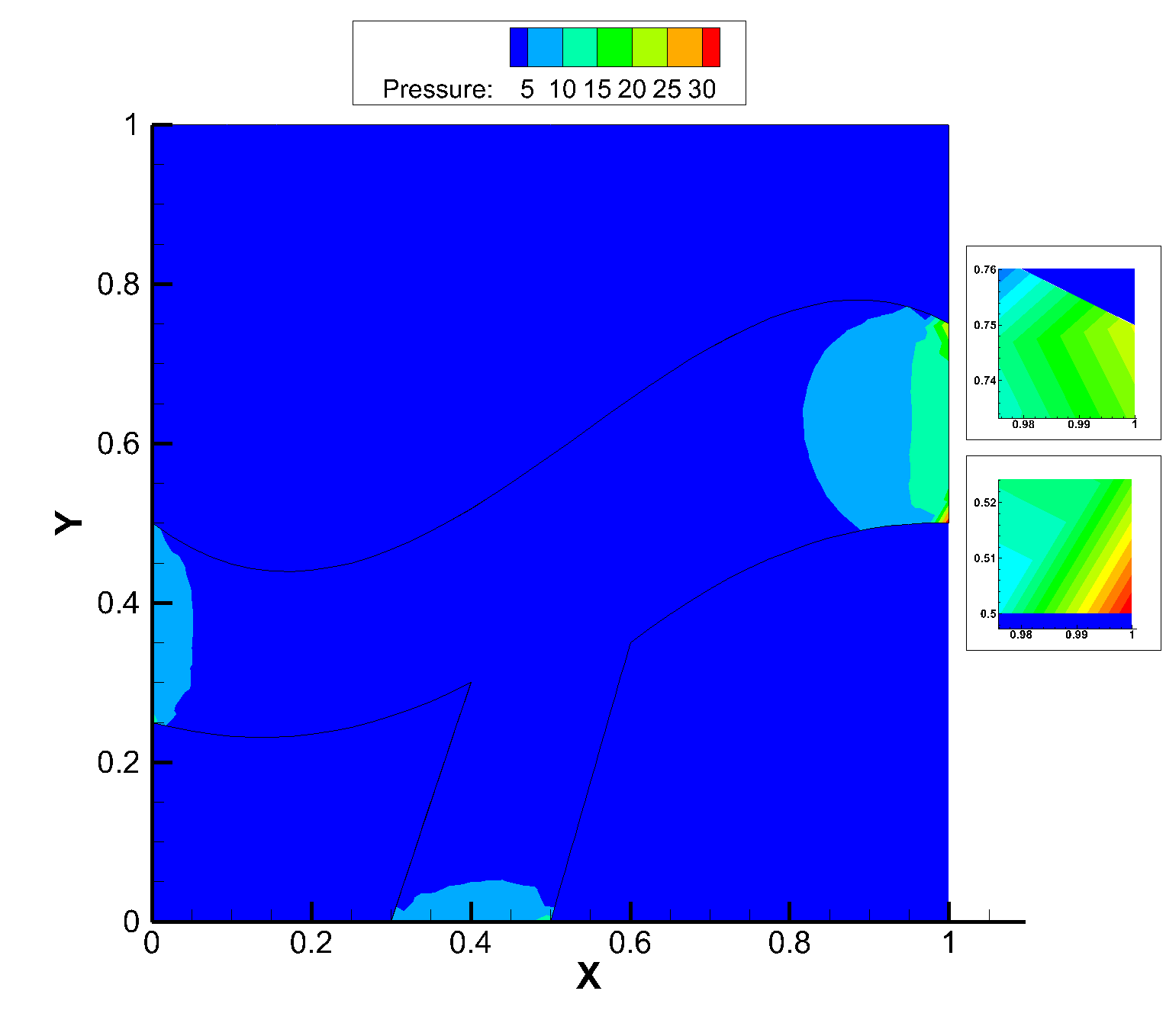}}
  \end{minipage}
  \caption{ Evolutions of velocity and pressure for Algorithm 1 ($\omega_1=0.5,~\omega_2=1,~k=10^{-2}$).}\label{k=-2}
\end{figure}

\begin{figure}[htbp]
  \centering
  \begin{minipage}[t]{0.49\linewidth}
  \centering
  \subfloat[][Velocity]{\label{equalVk-4}
  \includegraphics[width=7.1cm]{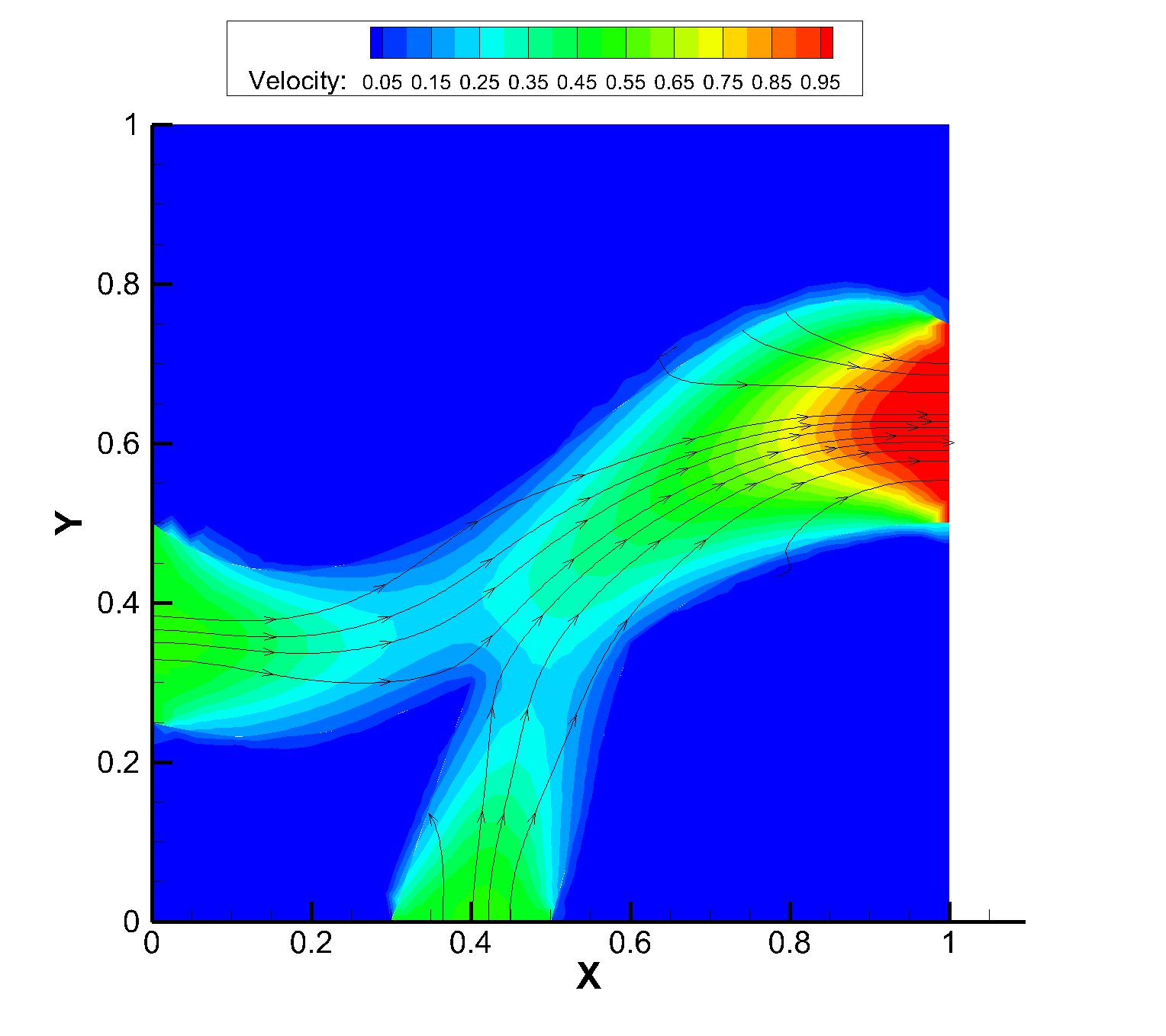}}
  \end{minipage}
  \begin{minipage}[t]{0.49\linewidth}
  \centering
  \subfloat[][Pressure]{\label{equalPk-4}
  \includegraphics[width=7.1cm]{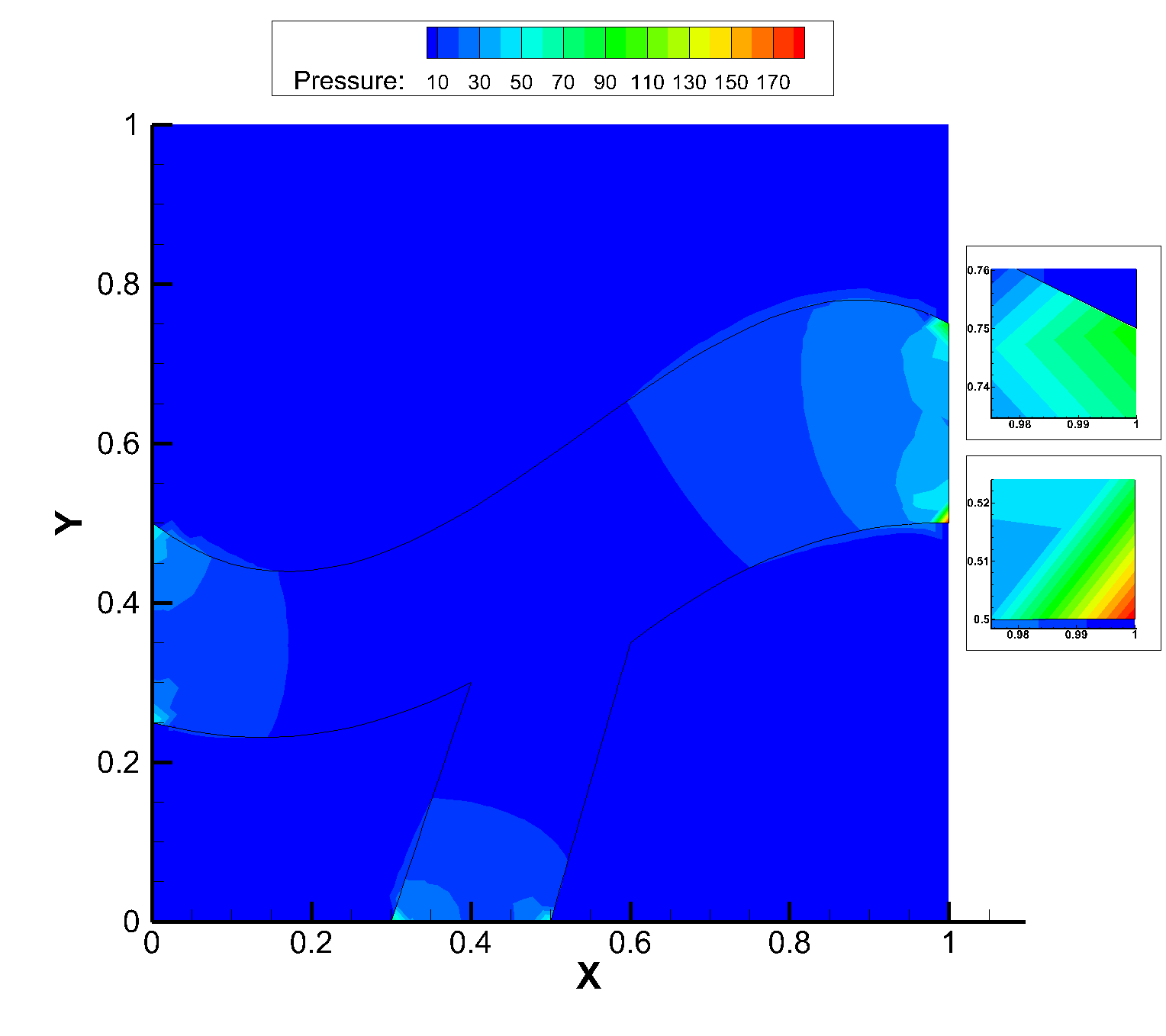}}
  \end{minipage}
  \caption{Evolutions of velocity and pressure for  Algorithm 1 ($\omega_1=0.5,~\omega_2=1,~k=10^{-4}$).}\label{k=-4}
\end{figure}

We present the streamlines of the velocity and the pressure in Figures \ref{k=1}-\ref{k=-4} by using different hydraulic conductivity $k=1, 10^{-2}, 10^{-4}$ respectively. Please note that we take the total inflow rate to be equal to the total outflow rate in above simulations. It can be easily observed that the fluid velocity in the porous media region decreases as the hydraulic conductivity decreases, which is consistent with the simulations in \cite{jiang2021SAV,Qiu2020domain,li2017stabilized}.

\begin{figure}[htbp]
  \centering
  \begin{minipage}[t]{0.49\linewidth}
  \centering
  \subfloat[][The total inflow rate is higher than the total outflow\\ rate ($\omega_1=1,~\omega_2=1$,~$k=1$).]{\label{dayuVk1}
  \includegraphics[width=7.1cm]{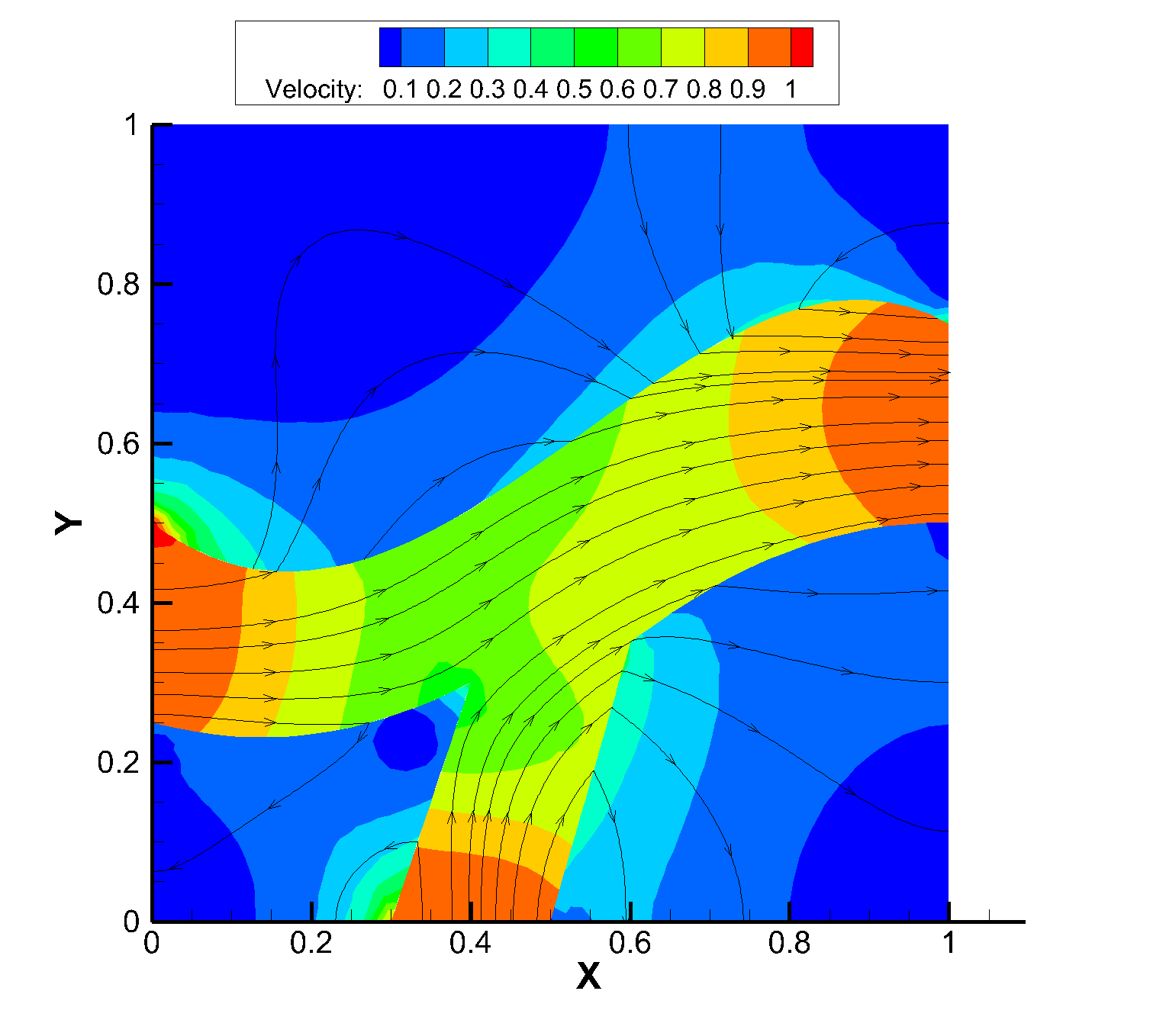}}
  \end{minipage}
  \begin{minipage}[t]{0.49\linewidth}
  \centering
  \subfloat[][The total inflow rate is lower than the total outflow \\ rate ($\omega_1=0.5,~\omega_2=2,~k=1$).]{\label{xiaoyuVk1}
  \includegraphics[width=7.1cm]{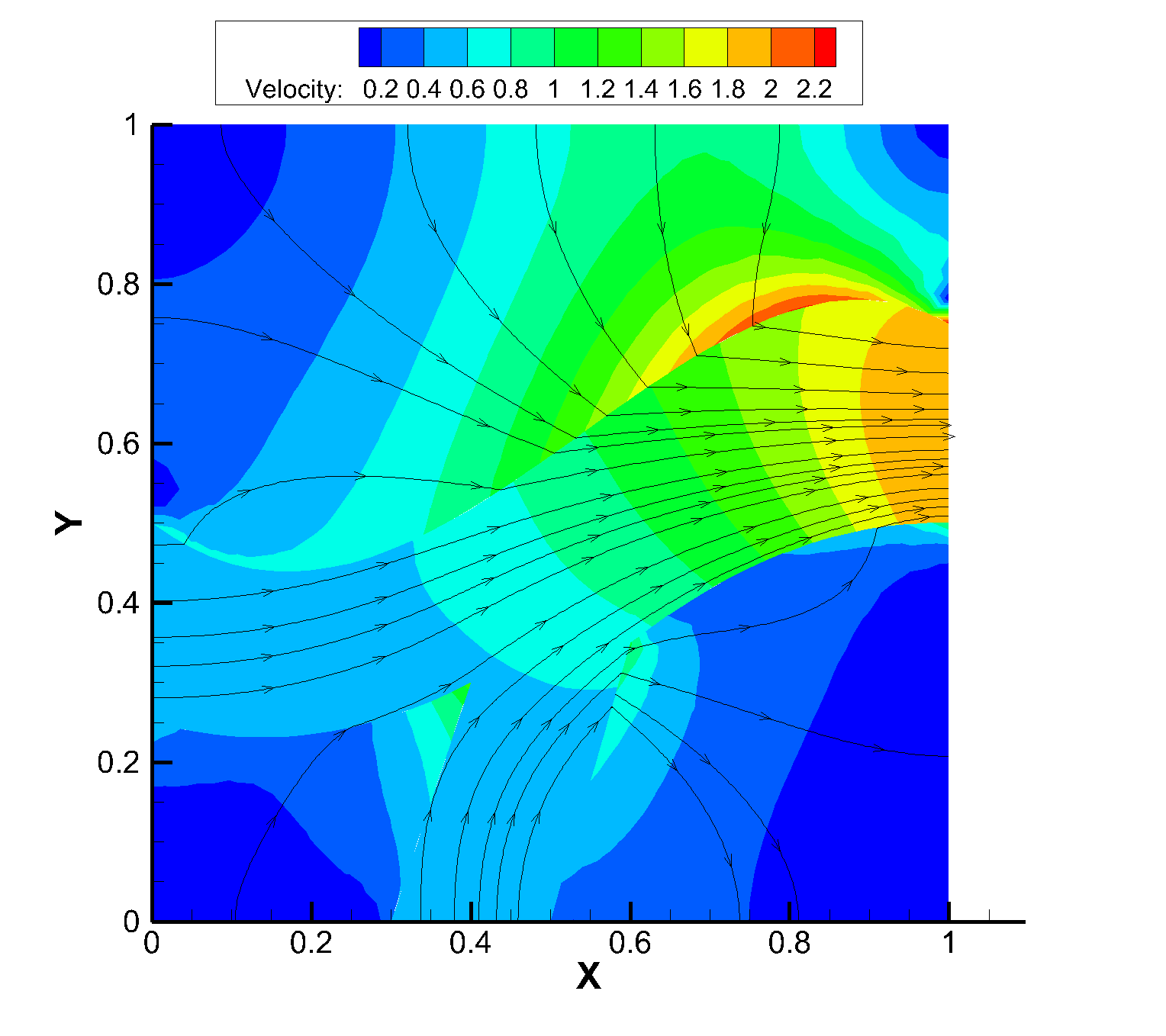}}
  \end{minipage}
  \caption{Evolutions of velocity when inputs and outputs are unequal.}\label{in out}
\end{figure}

Next, we will fix the hydraulic conductivity with $k=1$ and observe how the flow changes when the inflow velocity is not equal to the outflow velocity. As we can see in figure \ref{in out}, more flow will occur in the Y-shape region of the porous media when the total inflow rate is greater than the total outflow rate. The opposite is true if the total inflow rate is less than the total outflow rate, which is consistent with \cite{Qiu2020domain}.

\subsection{Filtration processes}
In this example, we investigate the impact of various hydraulic conductivities on fluid flow within porous media, which can be regarded as a simple simulation of a filter. In particular,
we will test the effect of the location and length of the low hydraulic conductivity regions on fluid flow by using Algorithm 1. Consider the computational domain as $\Omega_f=[0,2]\times [1.5,2]$, $\Omega_p=[0,2]\times [0,1.5]$ and $\Gamma=[0,2]\times \{1.5\}$. Set $T=0.5$, $\Delta t=0.01$, $h=\frac{1}{32}$, and other parameters all be 1. The hydraulic conductivity is $10^{-6}$ at two small rectangular portions of the porous media region and 1 elsewhere. The boundary conditions are given as
\begin{equation*}
\left\{
\begin{aligned}
&u=(0,x(x-2)),\quad{\rm on}~[0,2]\times\{2\},\\
&u=(0,0),\qquad\qquad{\rm on}~\{0\}\times [1.5,2]\cup \{2\}\times [1.5,2],\\
&\nabla\phi \cdot n_p=0,\qquad\quad{\rm on}~\{0\}\times[0,1.5] \cup \{2\}\times [0,1.5],\\
&\phi=0,\qquad\qquad\quad~~{\rm on}~[0,2]\times\{0\}.
\end{aligned}\right.
\end{equation*}

\begin{figure}[ht]
\centering
  \begin{minipage}[t]{0.33\linewidth}
    \centering
    \subfloat[][]{\label{11}
    \includegraphics[width=5.2cm]{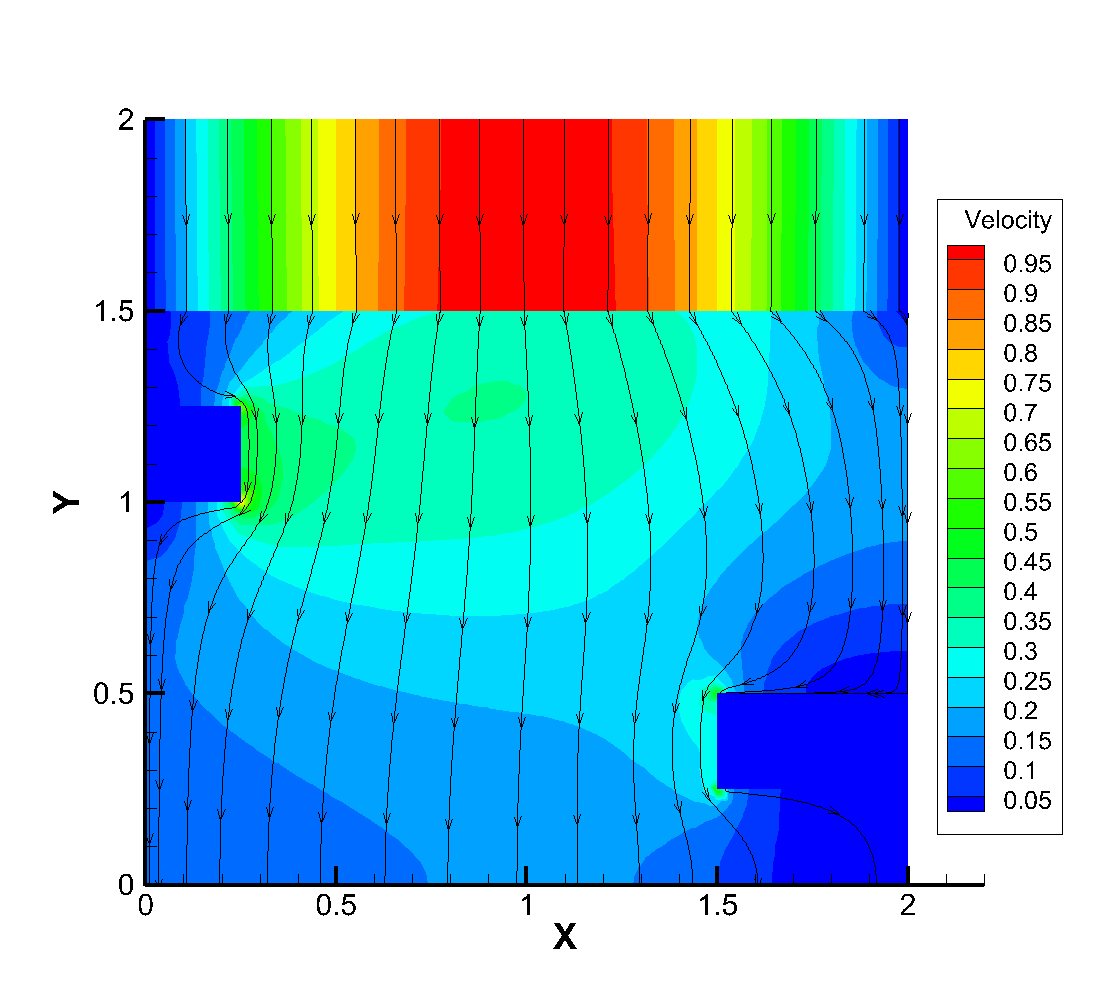}}
  \end{minipage}%
  \begin{minipage}[t]{0.33\linewidth}
    \centering
    \subfloat[][]{\label{12}
    \includegraphics[width=5.2cm]{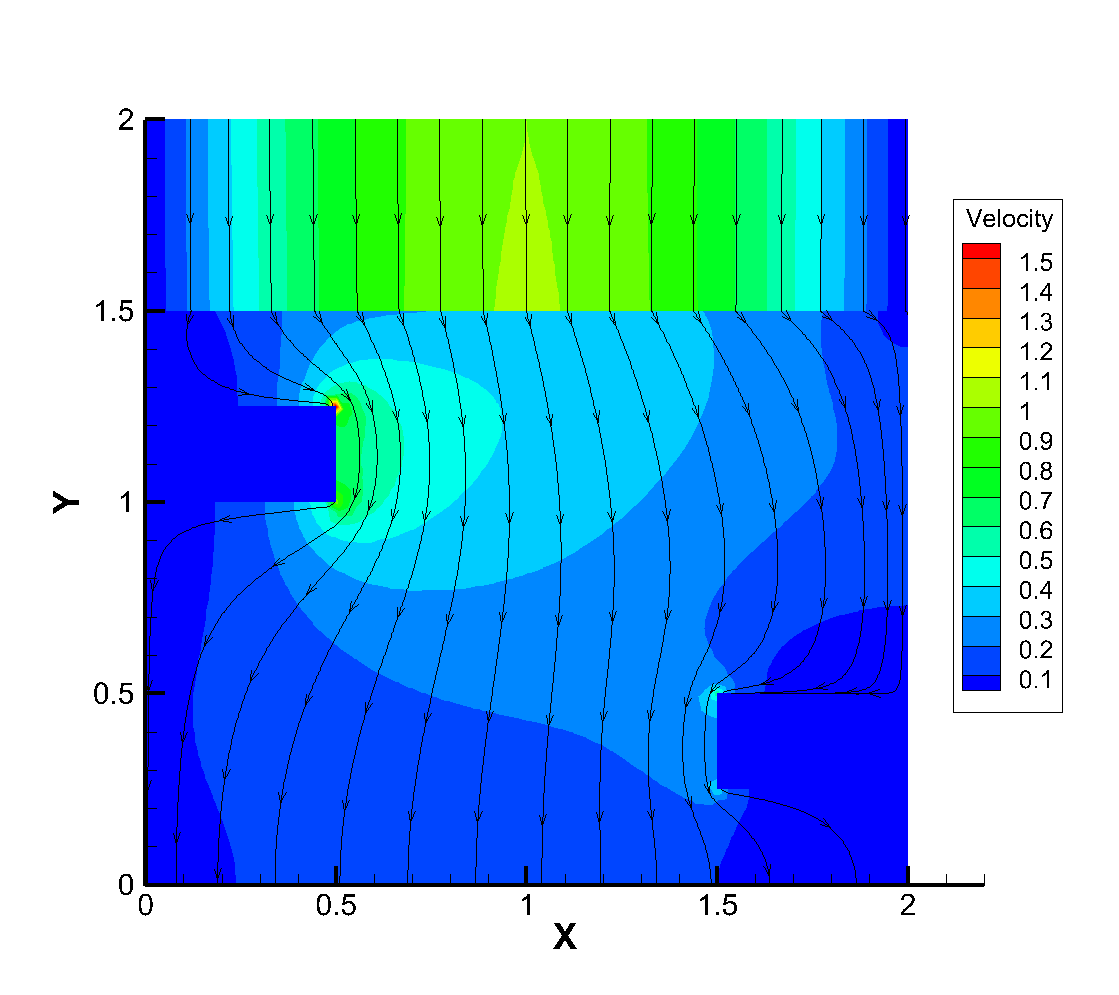}}
  \end{minipage}
  \begin{minipage}[t]{0.33\linewidth}
  \subfloat[][]{\label{13}
    \includegraphics[width=5.2cm]{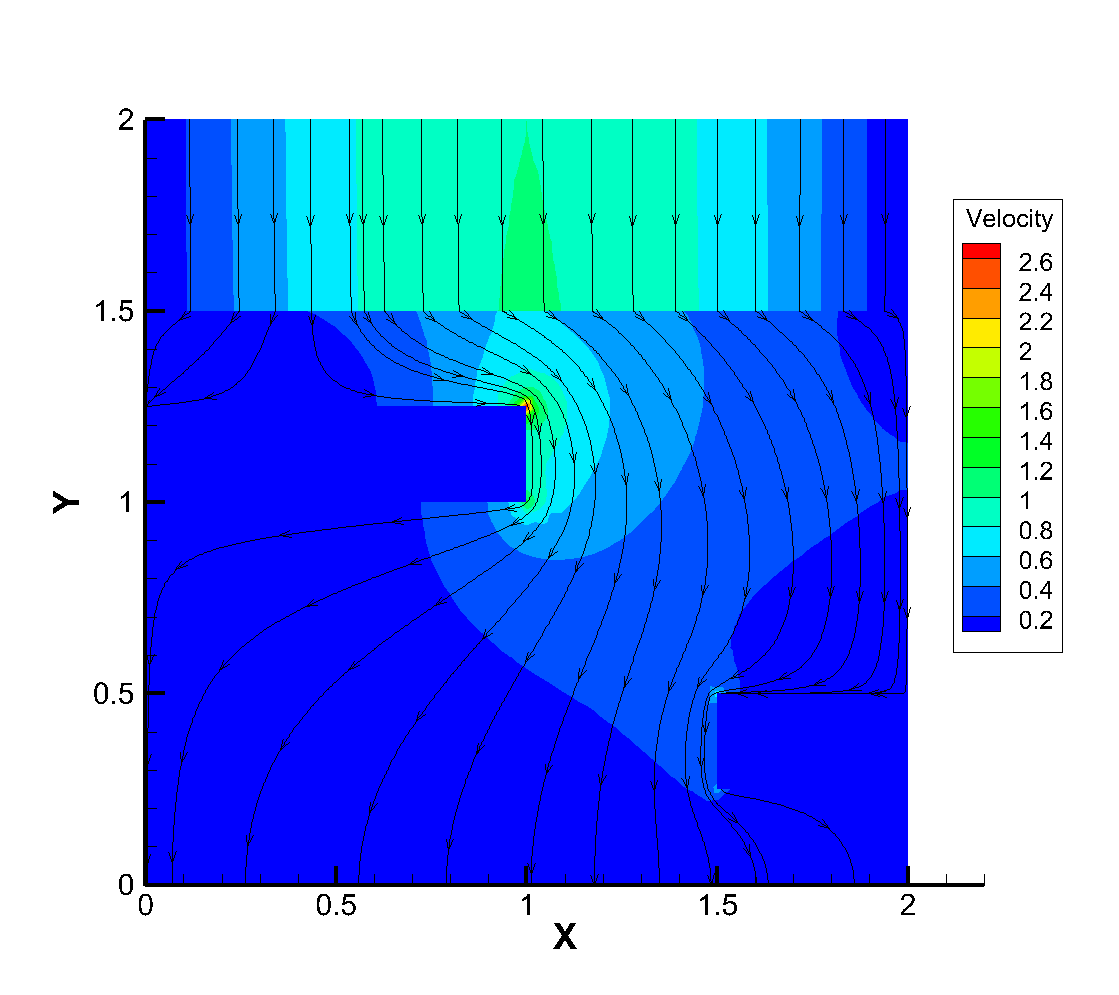}}
  \end{minipage}\\
   \begin{minipage}[t]{0.33\linewidth}
    \centering
    \subfloat[][]{\label{21}
    \includegraphics[width=5.2cm]{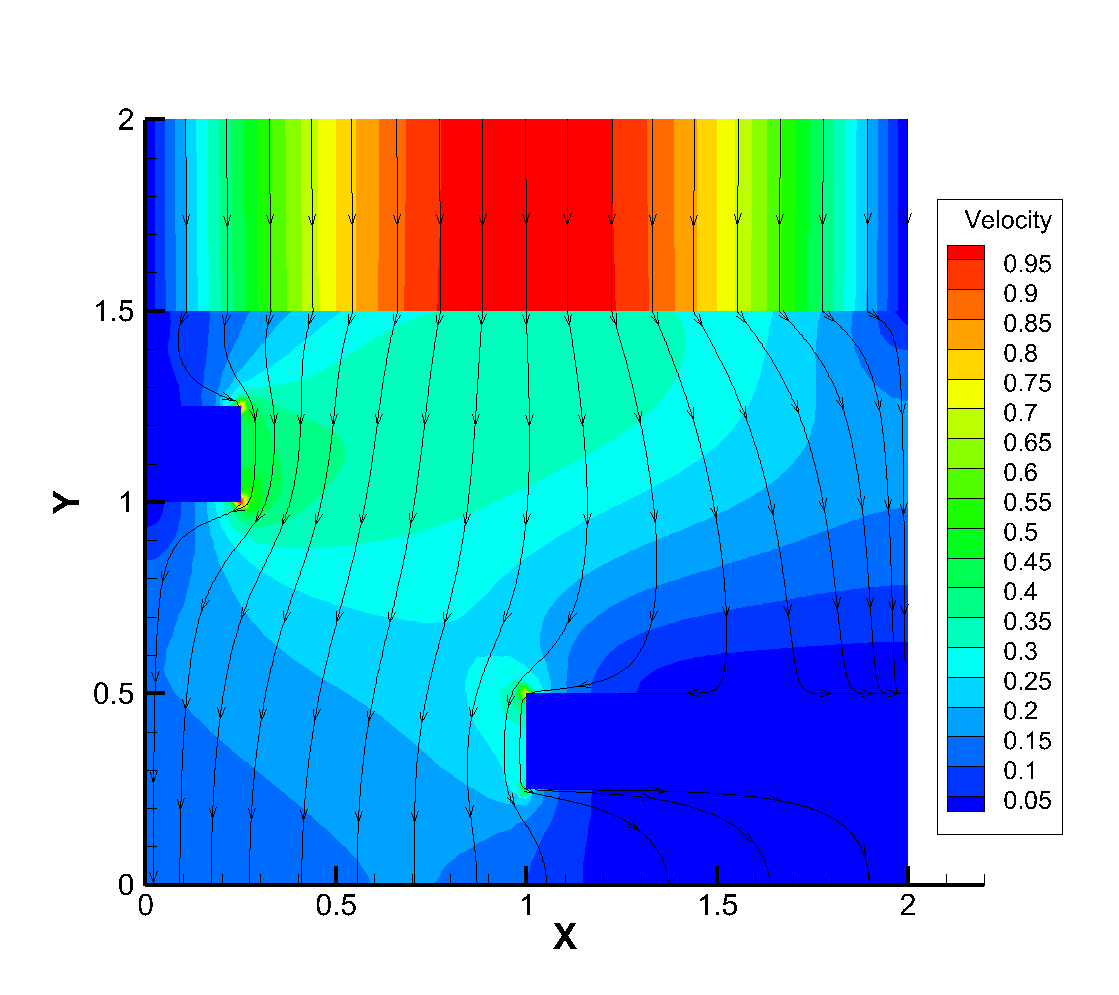}}
  \end{minipage}%
  \begin{minipage}[t]{0.33\linewidth}
    \centering
    \subfloat[][]{\label{22}
    \includegraphics[width=5.2cm]{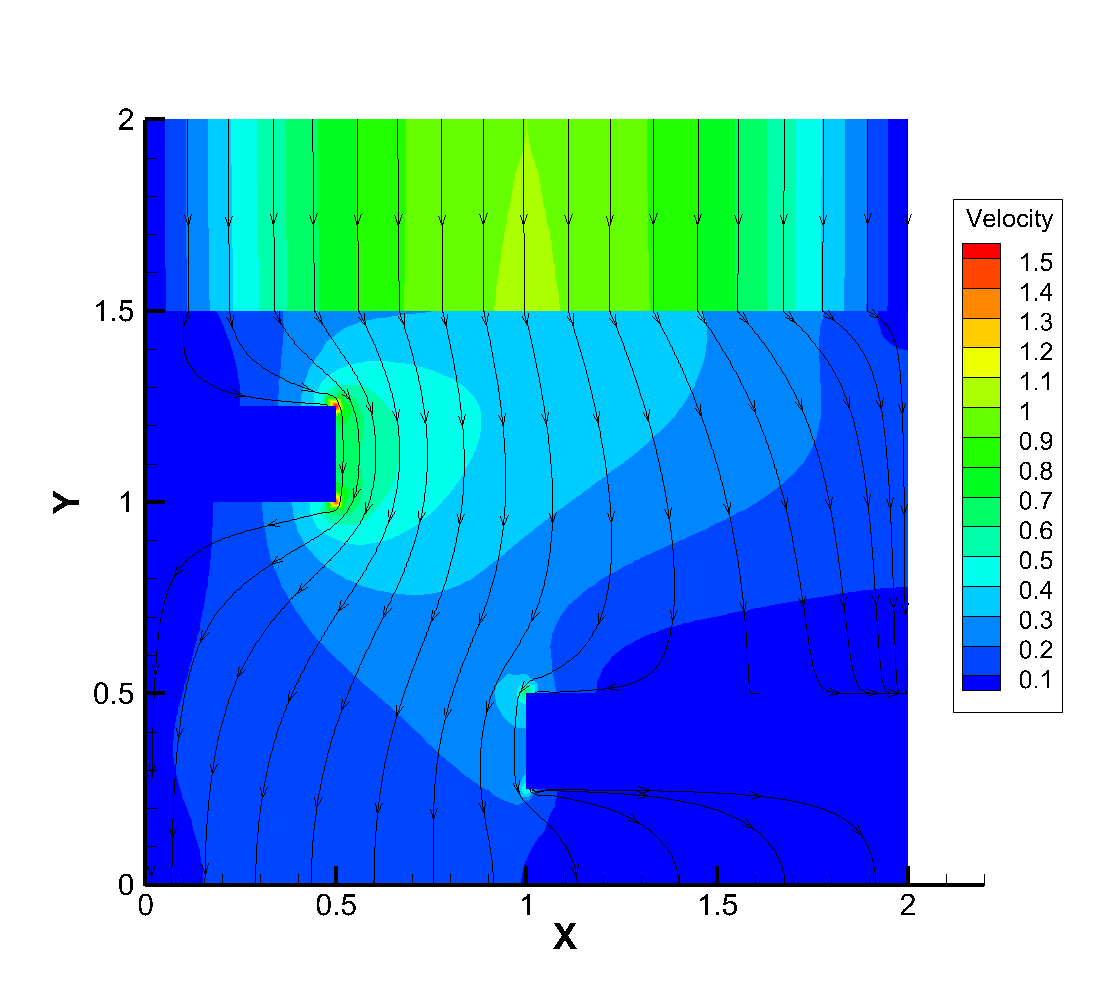}}
  \end{minipage}
  \begin{minipage}[t]{0.33\linewidth}
  \subfloat[][]{\label{23}
    \includegraphics[width=5.2cm]{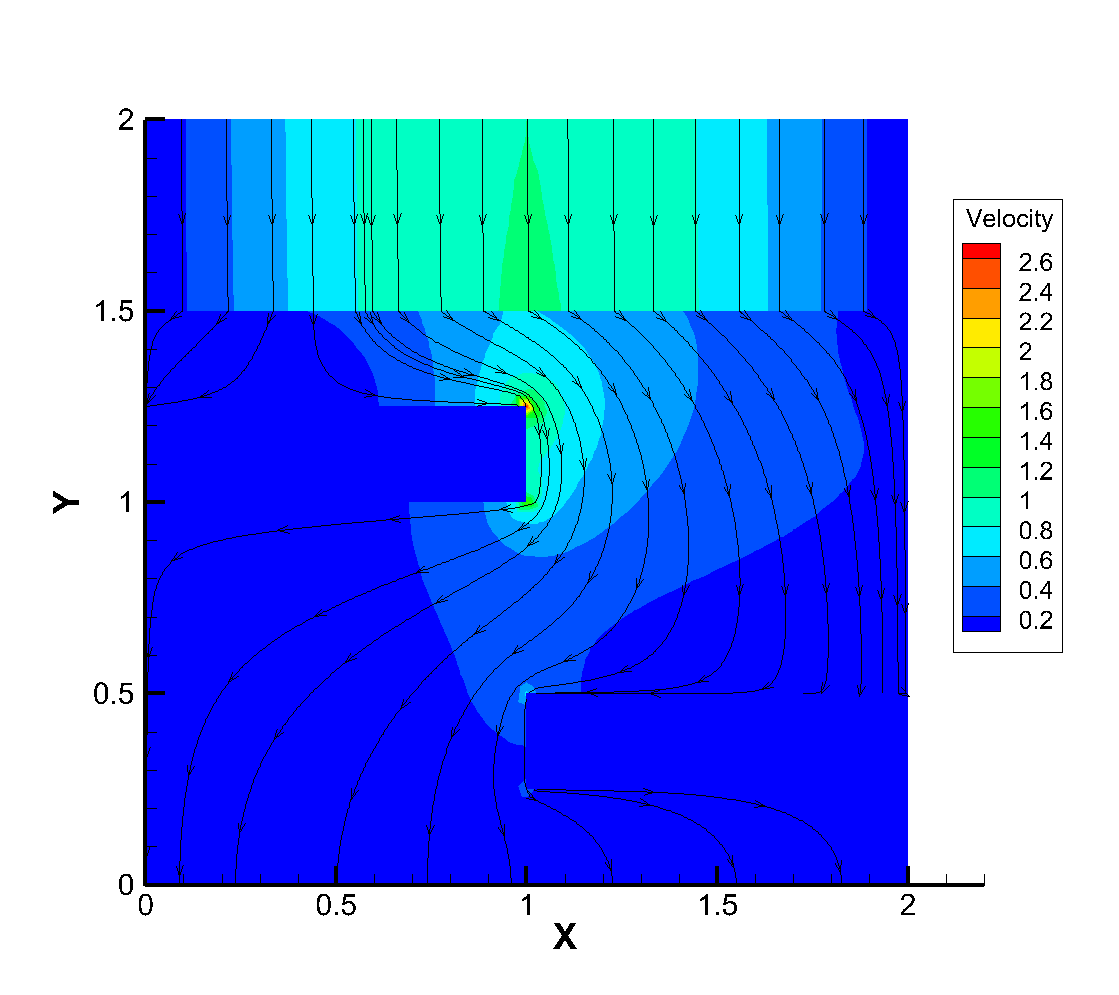}}
  \end{minipage}
  \caption{Comparison of velocity fields and streamlines for different filter settings.}
    \label{filtration comaparison}
\end{figure}

As we can see in Figure \ref{filtration comaparison}, the fluid first flows from $\Omega_f$ to $\Omega_p$ at a given rate, and as it passes through the first low hydraulic conductivity region, it follows its shape and accelerates at the corners. A similar phenomenon also occurs in the second low hydraulic conductivity block. In the case of the horizontal comparison, by fixing the width of the second region and changing the size of the first rectangular region, it is evident that acceleration becomes more pronounced as the width of the low hydraulic conductivity region increases, resulting in higher velocities. We can see a similar appearance by taking the comparison column by column. This phenomenon is similar to the narrow tube effect in airflow.

\subsection{Cavity flow} 
\begin{figure}[htbp]
\centering
  \begin{minipage}[t]{0.49\linewidth}
    \centering
    \subfloat[][Algorithm 1]{\label{cavitySAV}
    \includegraphics[width=7.1cm]{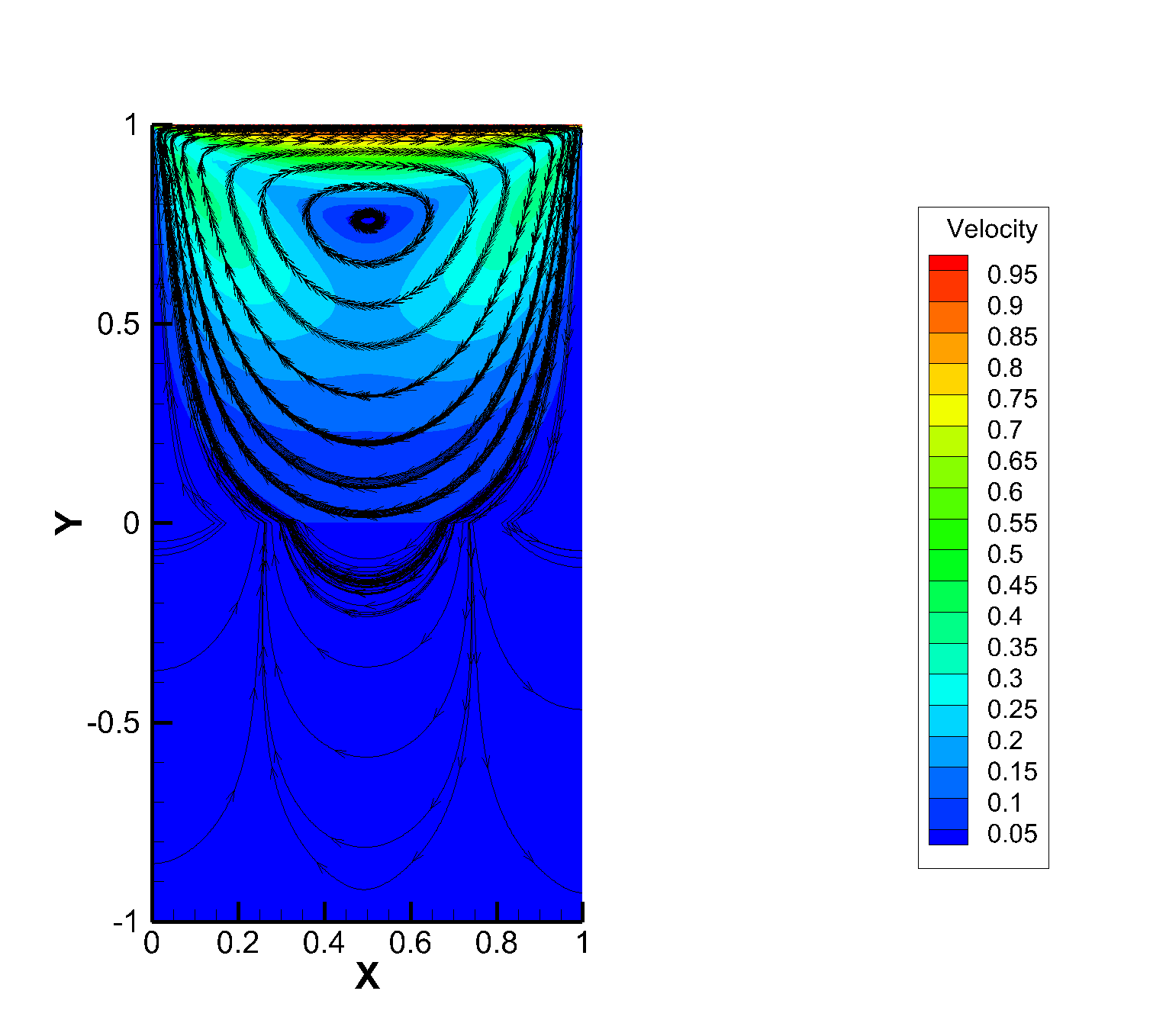}}
  \end{minipage}%
  \begin{minipage}[t]{0.49\linewidth}
    \centering
    \subfloat[][NDFEM]{\label{cavityNDFEM}
    \includegraphics[width=7.1cm]{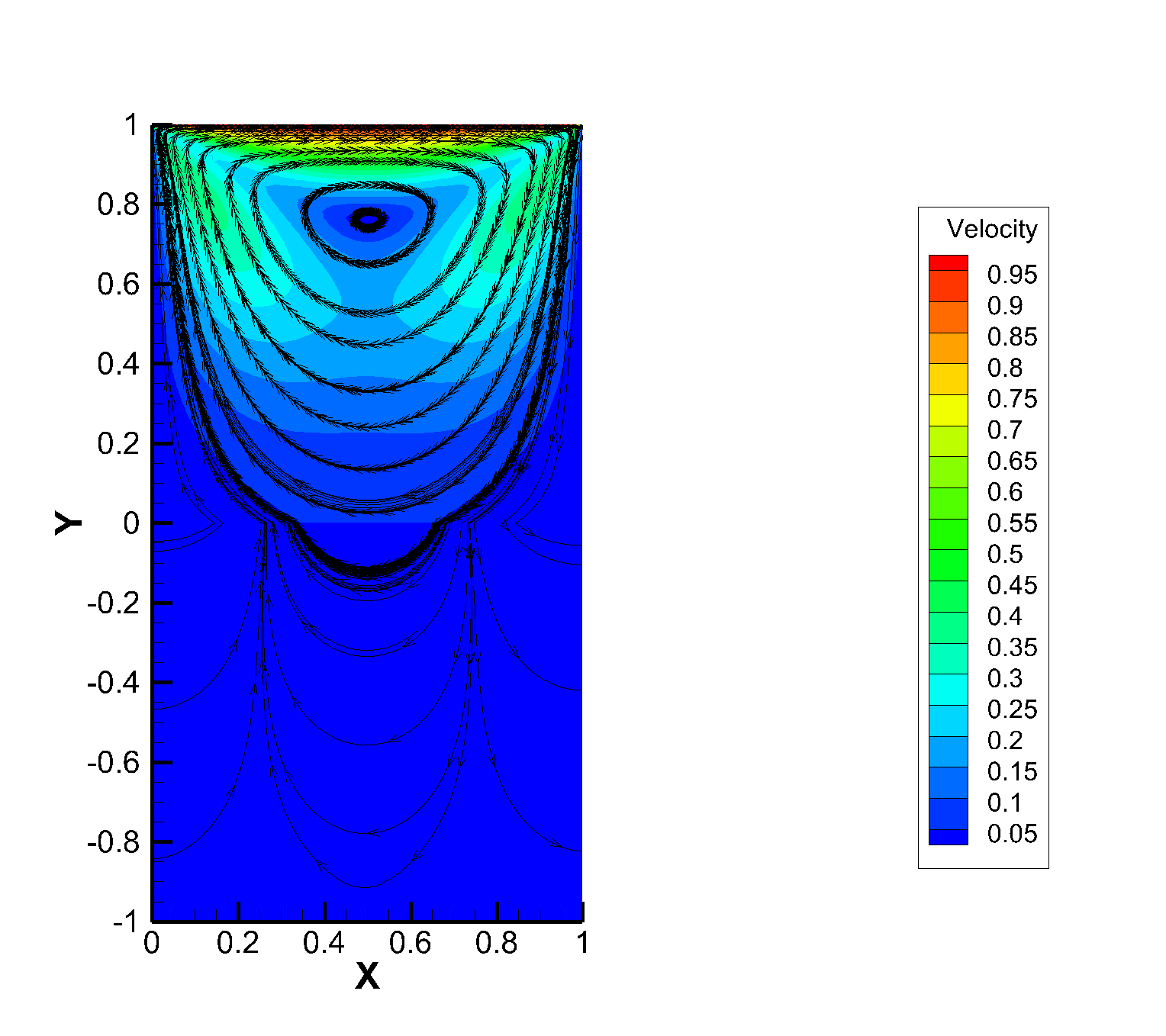}}
  \end{minipage}
    \caption{Comparisons of velocity field and streamlines of cavity flow.}
    \label{cavity comaparison}
\end{figure}

\begin{figure}[htbp]
\centering
  \begin{minipage}[t]{0.49\linewidth}
    \centering
    \subfloat[][Velocity component U1]{\label{cavitySAVDEMU1}
    \includegraphics[width=7.1cm]{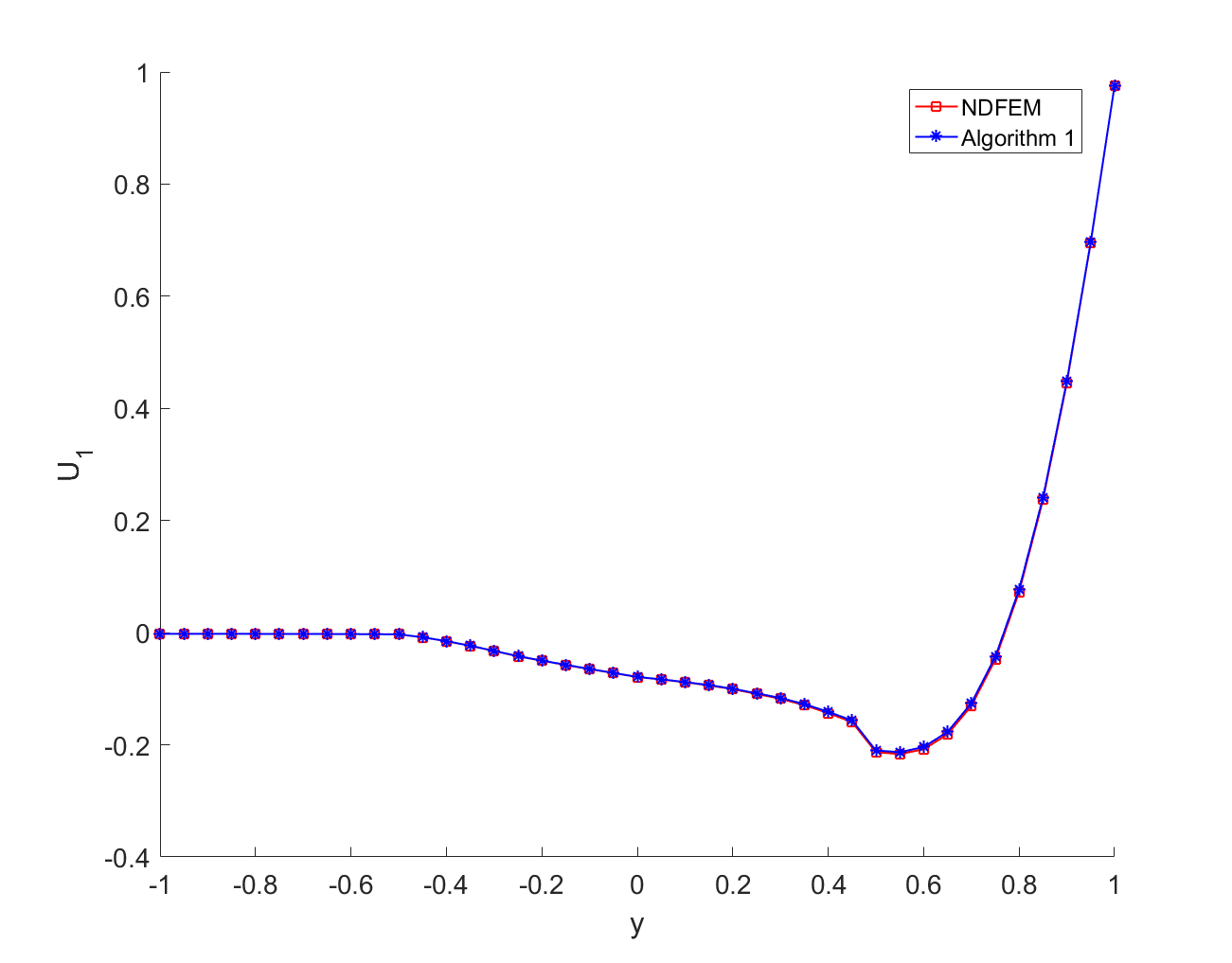}}
  \end{minipage}%
  \begin{minipage}[t]{0.49\linewidth}
    \centering
    \subfloat[][Velocity component U2]{\label{cavitySAVDEMU2}
    \includegraphics[width=7.1cm]{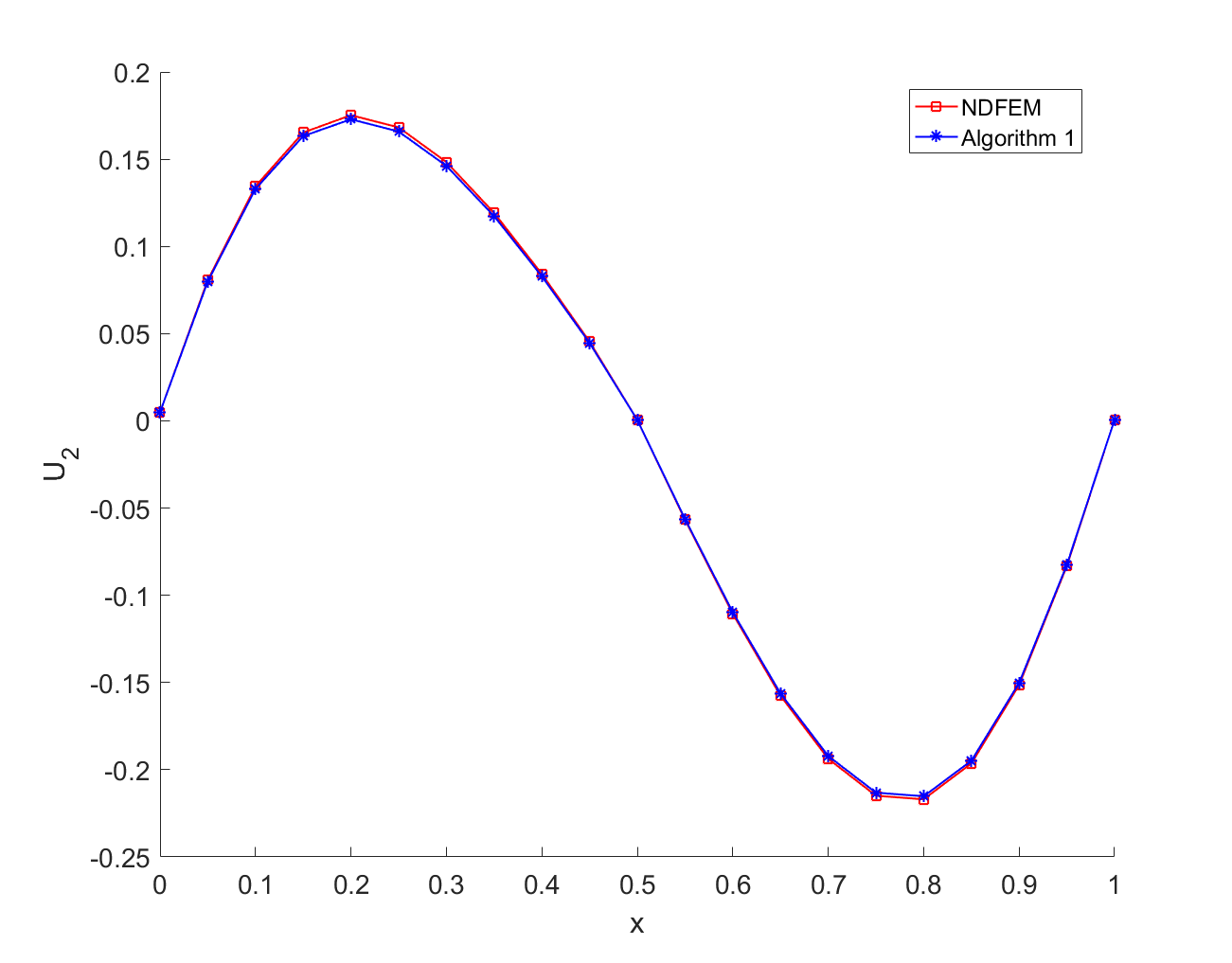}}
  \end{minipage}
    \caption{Numerical solutions of velocity by using  Algorithm 1 and NDFEM.}
    \label{cavity comaparison lines}
\end{figure}

In this example we will give the cavity flow test. The computational regions are defined as $\Omega_f=[0,1]\times [0,1]$, $\Omega_p=[0,1]\times [-1,0]$ and $\Gamma=[0,1]\times\{0\}$. Take $T=0.5$, $\Delta t=0.01$, $h=\frac{1}{64}$, all other parameters are 1. Given the following boundary conditions:
\begin{equation*}
\left\{
\begin{aligned}
&u=(1,0),\quad {\rm on}~[0,1]\times \{1\},\\
&u=(0,0),\quad {\rm on}~\{0\}\times [0,1]\cup \{1\}\times [0,1],\\
&\phi=0,\qquad~~ {\rm on}~\partial\Omega_p\backslash \Gamma.
\end{aligned}\right.
\end{equation*}
Define the global velocity $U=(U_1,U_2)=(u_1-k\frac{\partial \phi}{\partial x},u_2-k\frac{\partial \phi}{\partial y})$. We first present the velocity fields and streamlines of cavity flow by using Algorithm 1 and NDFEM in Figure \ref{cavity comaparison} where NDFEM denotes that using backward Euler scheme for the temporal discretization and Newton iteration for the trilinear term with decoupling finite element method, and then show that the velocity component $U_1$ obtained by both algorithms follows the same trend as $y$ increases from $-1$ to $1$ when $x$ is fixed at $0.5$ and similarly, the velocity component $U_2$ obtained by both algorithms follows a similar trend as $x$ increases from 0 to 1 with fixed $y = 0.5$ in  Figure \ref{cavity comaparison lines}, which verifies the accuracy and efficiency of our proposed method.

\section{Conclusion}
We proposed the efficient first- and second-order fully discrete schemes for the unsteady Navier-Stokes-Darcy model by using the SAV approach based on the finite element method. The nonlinear terms are explicitly discretized which means that only the linear differential equations need to be solved at each time step. We derived the unconditional stability of the first- and second-order numerical methods by modifying the trilinear term of the coupled system equipped with the Lions interface condition to preserve that $a_N(u,u,u)=0$. More importantly, we also established rigorous error estimates for the velocity and hydraulic head of the first-order scheme without any time step constraint by using the built-in stability result.

\section*{Appendix \uppercase\expandafter{\romannumeral1}\quad The proof of (\ref{err 2dt})}
\setcounter{equation}{0}
\renewcommand{\theequation}{\uppercase\expandafter{\romannumeral1}.\arabic{equation}}

\begin{flalign}\label{R1}
&|\bm{R_1}|=|a_N(u^{n+1}-u_h^n,u^{n+1},\theta_u^{n+1})|
\leq C\|u^{n+1}-u_h^n\|_{0,f}\|u^{n+1}\|_{2,f}\|\nabla\theta_u^{n+1}\|_{0,f},\nonumber\\
&\quad\quad\leq \varepsilon_1\nu\|\nabla\theta_u^{n+1}\|_{0,f}^2
+C\Delta t\int_{t^n}^{t^{n+1}}\|\partial _tu\|_{0,f}^2dt
+C\|\theta_u^n\|_{0,f}^2+C\|\sigma_u^{n}\|_{0,f}^2, &
\end{flalign}
\begin{flalign}\label{R2}
&|\bm{R_2}|=|a_N(u_h^n,u^{n+1}-u_h^n,\theta_u^{n+1})|\nonumber\\
&\qquad=|a_N(u_h^n,u^{n+1}-u^n,\theta_u^{n+1})+a_N(u^n-e_u^n,e_u^n,\theta_u^{n+1})|\nonumber\\
&\qquad=|a_N(u_h^n,u^{n+1}-u^n,\theta_u^{n+1})+a_N(u^n,\theta_u^n,\theta_u^{n+1})
+a_N(u^n,\sigma_u^{n},\theta_u^{n+1})
-a_N(\theta_u^n,\theta_u^n,\theta_u^{n+1})\nonumber\\
&\qquad\quad-a_N(\theta_u^n,\sigma_u^n,\theta_u^{n+1})
-a_N(\sigma_u^n,\theta_u^n,\theta_u^{n+1})-a_N(\sigma_u^n,\sigma_u^n,\theta_u^{n+1})|\nonumber\\
&\qquad\leq C\|u_h^n\|_{0,f}\|u^{n+1}-u^n\|_{2,f}\|\nabla\theta_u^{n+1}\|_{0,f}
+C\|u^n\|_{1,f}\|\theta_u^n\|_{0,f}^{1/2}\|\nabla\theta_u^n\|_{0,f}^{1/2}
\|\theta_u^{n+1}\|_{0,f}^{1/2}
\|\nabla \theta_u^{n+1}\|_{0,f}^{1/2}\nonumber\\
&\qquad\quad+C\|u^n\|_{1,f}\|\nabla\sigma_u^n\|_{0,f}\|\nabla \theta_u^{n+1}\|_{0,f}
+C\|\theta_u^n\|_{0,f}^{1/2}\|\nabla\theta_u^n\|_{0,f}^{1/2}
\|\theta_u^{n}\|_{0,f}^{1/2}\|\nabla\theta_u^n\|_{0,f}^{1/2}
\|\nabla\theta_u^{n+1}\|_{0,f}\nonumber\\
&\qquad
\quad+C\|\nabla\sigma_u^n\|_{0,f}\|\theta_u^n\|_{0,f}^{1/2}
\|\nabla\theta_u^n\|_{0,f}^{1/2}\|\theta_u^{n+1}\|_{0,f}^{1/2}
\|\nabla\theta_u^{n+1}\|_{0,f}^{1/2}
+C\|\nabla\sigma_u^n\|_{0,f}\|\nabla\sigma_u^n\|_{0,f}\|\nabla\theta_u^{n+1}\|_{0,f}\nonumber\\
&\qquad
\leq C\Delta t\int_{t^n}^{t^{n+1}}\|\partial _t u\|_{2,f}^2dt
+\varepsilon_1\nu\|\nabla \theta_u^{n+1}\|_{0,f}^2
+C\|\theta_u^n\|_{0,f}\|\nabla\theta_u^n\|_{0,f}
+C\|\theta_u^{n+1}\|_{0,f}\|\nabla\theta_u^{n+1}\|_{0,f}
\nonumber\\
&\qquad\quad+C\|\nabla\sigma_u^n\|_{0,f}^2
+C\|\theta_u^n\|_{0,f}^2\|\nabla\theta_u^n\|_{0,f}^2
+C\|\nabla\sigma_u^n\|_{0,f}^2\|\theta_u^n\|_{0,f}\|\nabla\theta_u^n\|_{0,f}
+C\|\theta_u^{n+1}\|_{0,f}\|\nabla\theta_u^{n+1}\|_{0,f}\nonumber\\
&\qquad\quad+C\|\nabla\sigma_u^n\|_{0,f}^4
\nonumber\\
&\qquad
\leq C\Delta t\int_{t^n}^{t^{n+1}}\|\partial_tu\|_{2,f}^2dt
+\varepsilon_1\nu\|\nabla\theta_u^{n+1}\|_{0,f}^2+C\|\theta_u^n\|_{0,f}^2
+\varepsilon_1\nu\|\nabla\theta_u^{n}\|_{0,f}^2+C\|\theta_u^{n+1}\|_{0,f}^2
\nonumber\\
&\qquad\quad+C\|\nabla\sigma_u^n\|_{0,f}^2
+C\|\theta_u^n\|_{0,f}^2\|\nabla\theta_u^n\|_{0,f}^2
+
{C\|\nabla\sigma_u^n\|_{0,f}^4}, &
\end{flalign}
\begin{flalign}\label{R3}
&|\bm{R_3}|=|c_\Gamma(\theta_u^{n+1},\phi^{n+1}-\phi_h^n)|\nonumber\\
&\qquad\leq C\|\nabla\theta_u^{n+1}\|_{0,f}\|\phi^{n+1}-\phi_h^n\|_{0,p}^{1/2}
\|\nabla(\phi^{n+1}-\phi_h^n)\|_{0,p}^{1/2}\nonumber\\
&\qquad\leq C\|\phi^{n+1}-\phi_h^n\|_{0,p}\|\nabla(\phi^{n+1}-\phi_h^n)\|_{0,p}
+\varepsilon_1\nu\|\nabla\theta_u^{n+1}\|_{0,f}^2\nonumber\\
&\qquad\leq C\|\phi^{n+1}-\phi_h^n\|_{0,p}^2
+\varepsilon \|\nabla(\phi^{n+1}-\phi_h^n)\|_{0,p}^2
+\varepsilon_1\nu\|\nabla\theta_u^{n+1}\|_{0,f}^2\nonumber\\
&\qquad\leq C\Delta t\int_{t^n}^{t^{n+1}}\|\partial_t \phi\|_{0,p}^2dt
+C\|\theta_\phi^n\|_{0,p}^2+C\|\sigma_\phi^n\|_{0,p}^2
+C\Delta t\int_{t^n}^{t^{n+1}}\|\partial_t \phi\|_{1,p}^2dt\nonumber\\
&\qquad\quad+\varepsilon_2 gk\|\nabla\theta_\phi^n\|_{0,p}^2
+C\|\nabla\sigma_\phi^n\|_{0,p}^2
+\varepsilon_1\nu\|\nabla\theta_u^{n+1}\|_{0,f}^2, &
\end{flalign}
\begin{flalign}\label{R4}
&|\bm{R_4}|=|c_\Gamma(u^{n+1}-u_h^n,\theta_\phi^{n+1})|\nonumber\\
&\qquad\leq C\|\nabla\theta_\phi^{n+1}\|_{0,p}\|u^{n+1}-u_h^n\|_{0,f}^{1/2}
\|\nabla(u^{n+1}-u_h^n)\|_{0,f}^{1/2}\nonumber\\
&\qquad\leq \varepsilon_2gk\|\nabla\theta_\phi^{n+1}\|_{0,p}^2
+C\|u^{n+1}-u_h^n\|_{0,f}^2
+\varepsilon\|\nabla(u^{n+1}-u_h^n)\|_{0,f}^2\nonumber\\
&\qquad\leq \varepsilon_2gk\|\nabla\theta_\phi^{n+1}\|_{0,p}^2
+C\Delta t\int_{t^n}^{t^{n+1}}\|\partial_t u\|_{0,f}^2dt
+C\|\theta_u^n\|_{0,f}^2+C\|\sigma_u^n\|_{0,f}^2
\nonumber\\
&\qquad\quad+C\Delta t\int_{t^n}^{t^{n+1}}\|\partial_t u\|_{1,f}^2dt
+\varepsilon_1\nu\|\nabla\theta_u^n\|_{0,f}^2
+C\|\nabla\sigma_u^n\|_{0,f}^2, &
\end{flalign}
\begin{flalign}\label{R5}
&\bm{R_5}=R_u^{n+1}(\theta_u^{n+1})=(\frac{u^{n+1}-u^n}{\Delta t}-\partial_tu(t^{n+1}),\theta_u^{n+1})_f
\leq C\Delta t\int_{t^n}^{t^{n+1}}\|\partial_{tt}u\|_{0,f}^2dt
+\varepsilon_1\nu\|\nabla\theta_u^{n+1}\|_{0,f}^2,&
\end{flalign}
\begin{flalign}\label{R6}
&\bm{R_6}=R_\phi^{n+1}(\theta_\phi^{n+1})
=gS_0(\frac{\phi^{n+1}-\phi^n}{\Delta t}-\partial_t\phi(t^{n+1}),\theta_\phi^{n+1})_p\nonumber\\
&\qquad\qquad\qquad\qquad
\leq C\Delta t\int_{t^n}^{t^{n+1}}\|\partial_{tt}\phi\|_{0,p}^2dt
+\varepsilon_2gk\|\nabla\theta_\phi^{n+1}\|_{0,p}^2,&
\end{flalign}
%
%
%
\begin{flalign}\label{R7}
&|\bm{R_7}|=|(\frac{\sigma_u^{n+1}-\sigma_u^n}{\Delta t},\theta_u^{n+1})_f
+gS_0(\frac{\sigma_\phi^{n+1}-\sigma_\phi^n}{\Delta t},\theta_\phi^{n+1})_p|\nonumber\\
&\qquad\leq
\frac{C}{\Delta t}\int_{t^n}^{t^{n+1}}\|\partial_t\sigma_u\|_{0,f}^2dt
+\varepsilon_1\nu\|\nabla\theta_u^{n+1}\|_{0,f}^2
+\frac{C}{\Delta t}\int_{t^n}^{t^{n+1}}\|\partial_t\sigma_\phi\|_{0,p}^2dt
+\varepsilon_2gk\|\nabla\theta_\phi^{n+1}\|_{0,p}^2,&
\end{flalign}

%
\begin{flalign}\label{R8}
&|\bm{R_{8}}|
=\left|\eta\sum_{i=1}^{d-1}\int_\Gamma(\sigma_u^{n+1}\cdot\tau_i)
(\theta_u^{n+1}\cdot\tau_i)ds\right|
\leq \eta\sum_{i=1}^{d-1}\left(C\|\nabla\sigma_u^{n+1}\|_{0,f}^2
+\varepsilon_3\int_\Gamma(\theta_u^{n+1}\cdot\tau_i)^2ds\right). &
\end{flalign}

Choosing $\varepsilon_3=\frac{1}{2}$, then combining \eqref{err ADD} with (\ref{R1})-(\ref{R8}) arrives at (\ref{err 2dt}).

\bibliographystyle{siamplain}
\bibliography{NSD_SAV}

\end{document}